\RequirePackage{amsthm}
\documentclass[pdflatex,sn-mathphys-num]{sn-jnl}

\usepackage{amsmath,amssymb,amsthm,mathtools,mathrsfs}
\usepackage{booktabs,multirow,array}
\usepackage{graphicx}
\usepackage{algorithm}
\usepackage[noend]{algpseudocode}
\usepackage{enumitem}
\usepackage{xcolor}
\hypersetup{colorlinks=true,linkcolor=blue!60!black,citecolor=blue!60!black,urlcolor=blue!60!black,hypertexnames=false}

\theoremstyle{thmstyleone}
\newtheorem{theorem}{Theorem}
\newtheorem{proposition}[theorem]{Proposition}
\newtheorem{lemma}[theorem]{Lemma}
\newtheorem{corollary}[theorem]{Corollary}
\theoremstyle{thmstyletwo}
\newtheorem{remark}{Remark}
\theoremstyle{thmstylethree}
\newtheorem{definition}{Definition}

\newlist{thmlist}{enumerate}{1}
\setlist[thmlist]{label=\textup{(\roman*)},ref=(\roman*),leftmargin=2.2em,itemsep=1pt,topsep=2pt}
\newcommand{\pitem}[1]{\leavevmode\par\noindent\textup{(#1)}\enspace\ignorespaces}

\newcommand{\tablesmall}{\fontsize{7bp}{8.5bp}\selectfont}
\setlist[itemize]{label=$\bullet$,leftmargin=1.5em,itemsep=1pt,topsep=3pt}
\newcommand{\solver}{LBP}

\providecommand{\R}{\mathbb{R}}
\providecommand{\Z}{\mathbb{Z}}
\providecommand{\E}{\mathbb{E}}
\providecommand{\Prob}{\mathbb{P}}
\providecommand{\Var}{\operatorname{Var}}
\providecommand{\Cov}{\operatorname{Cov}}
\providecommand{\Leh}{\mathscr{L}}
\providecommand{\dd}{\,\mathrm{d}}
\providecommand{\eqdef}{\mathrel{\vcentcolon=}}
\providecommand{\slack}{\operatorname{slack}}
\providecommand{\load}{\ell}
\providecommand{\calA}{\mathcal{A}}
\providecommand{\calC}{\mathcal{C}}

\providecommand{\zlp}{z_{\mathrm{LP}}}
\providecommand{\zrmp}{z_{\mathrm{RMP}}}
\providecommand{\LB}{\mathrm{LB}}
\providecommand{\UB}{\mathrm{UB}}
\providecommand{\rc}{\bar c}
\providecommand{\wv}{\mathbf{w}}
\providecommand{\dv}{\mathbf{d}}
\providecommand{\one}{\mathbf{1}}

\begin{document}

\title[Entropic Column Generation and Lehmer Branch-and-Price]{Entropic Column Generation and Lehmer Branch-and-Price for Bin Packing Problem}

\author*[1]{\fnm{Masoud} \sur{Ataei}}\email{masoud.ataei@utoronto.ca}
\affil*[1]{\orgdiv{Department of Mathematical and Computational Sciences}, \\ \orgname{University of Toronto}, \orgaddress{\country{Canada}}}

\abstract{We introduce Lehmer branch-and-price, an exact algorithm
for bin packing and cutting stock problems built on the Lehmer transform,
the ratio of Mellin transforms of item weights at consecutive orders. The
transform maps each bin, a finite multiset of item weights, to a smooth increasing curve,
and so gives the discrete problem a continuous representation on which the
tools of analysis, such as convexity, exponential families and smoothing,
can act. Read along this curve, the problem yields an explicit dual
solution, a monotone representation of the dominance order, a family of
constructive heuristics with first fit decreasing as its limit, and an
integer identity that tells when the linear relaxation can decide an
instance and when only a perfect packing can improve the best solution
known. Lifted from items to patterns, the transform gives entropic column
generation, a new form of column generation that bounds the relaxation
from both sides at every iteration without solving the pricing problem
exactly, maximizes a smoothed dual at an
explicit rate, and measures by an entropy how much the relaxation leaves
undecided. We show that the proposed algorithm outperforms the
state-of-the-art methods on the standard benchmark in both the number of
instances solved and the computational time.}


\maketitle

\clearpage
\section{Introduction}\label{sec:intro}

The bin packing problem asks for the fewest bins of a given capacity
into which a set of items of given weights can be packed; the cutting
stock problem is the same problem with the items grouped into types with
demands. Exact methods for both are organized around one integer program,
the pattern formulation of Gilmore and Gomory~\cite{gilmore1961}, whose
roots go back to Kantorovich and Zalgaller~\cite{uchoa2026}. It has one
variable for each pattern, a feasible filling of a bin, which counts how
many bins are filled that way, and one constraint per item type, which
requires its demand to be covered. Its linear relaxation is famously
strong: almost every instance has the integer round-up property (IRUP),
its optimum being the relaxation value rounded
up~\cite{scheithauer1995,kartak2015}, and the problem is polynomial for a
fixed number of item types~\cite{goemans2020}. Branch-and-price methods
solve the pattern formulation by column
generation~\cite{vance1994,vanderbeck1999,barnhart1998,belov2006},
pseudo-polynomial methods solve an arc-flow or related network
formulation with a general integer programming
solver~\cite{valerio1999,brandao2016,delima2022,delorme2020,delorme2026},
and the review of Delorme, Iori and Martello~\cite{delorme2016} fixed the
benchmark on which both lines are compared. In all of them a bin enters
only through its load and its item counts. It is a finite multiset
without smooth structure, and apart from dual-feasible functions, which
act on its item weights one at a time, the tools of analysis, convexity,
monotonicity, exponential families and smoothing, have little to act on.

This paper gives the bin that structure. The Lehmer transform of a bin,
the ratio of the Mellin transforms of its item weights at consecutive
orders~\cite{ataei2022,ataei2026,ataei_vinogradov2026}, turns the multiset
into an analytic curve that increases from its smallest to its largest
weight, and so gives a discrete problem a continuous representation. Read
at order one, the numerator of the transform is the load of the bin and
yields an explicit dual solution of the linear relaxation; read at higher
orders, it represents the dominance order monotonically and generates a
family of heuristics with first fit decreasing as its limit. With that dual, the
accounting identity of column generation becomes an identity in integers,
the slack identity, which tells when the relaxation can decide an instance
and when only a perfect packing can improve the best solution known. Lifted from items to patterns, the transform
becomes a probability distribution on the patterns, and column generation
acquires an entropic form: the free energy brackets the relaxation at
every iteration, a smoothed dual is maximized at an explicit rate, and its
duality gap is an entropy.

The hardest instances of the benchmark show what this structure buys. The
augmented IRUP and non-IRUP classes (AI and ANI) were built so that the
linear relaxation cannot tell the instances that attain its value from
those that need one more bin~\cite{delorme2016,caprara2015}. The leading
general-purpose codes, which rely on a commercial solver, leave some of
them open after an hour~\cite{wei2020,pessoa2020,delima2023,dasilva2025},
and da Silva et al.~\cite{dasilva2026} recently solved all of them by
algorithms specific to their construction. The slack identity explains
why: an improving solution must be a perfect packing, and the
relaxation, which values every full bin alike, cannot say which full bins
form one. Since the identity detects this on any instance, it tells an
exact method where combinatorial reasoning must take over.

On these foundations we build Lehmer branch-and-price (\solver). Every
bound it uses is certified, so every optimum it reports is proved. The
column generation of a node stops as soon as the node can be pruned or
must be branched on. Where only a perfect packing can improve the
incumbent, \solver{} switches to combinatorial reasoning: exchange
reductions fix bins that some perfect packing must contain, and backdoor
peeling, named after the backdoors of satisfiability and integer
programming~\cite{williams2003,fischetti2013}, holds back a few items so
that these reductions cascade. Unlike the leading general-purpose codes,
\solver{} uses no commercial solver. It solves every instance of the
benchmark within ten minutes, the augmented ones in a fraction of a
second on average.

Several ingredients are classical in other forms, and it is worth stating
precisely where the present work stands relative to them. The
monotonicity of the Lehmer means is classical~\cite{lehmer1971,bullen2003};
the characterization of dominance is the superadditive duality of
dual-feasible functions~\cite{clautiaux2010} applied to a pair of bins;
the slack identity is the waste argument made exact; the decision stop
refines the early termination of column generation by Lagrangian
bounds~\cite{lubbecke2005} and shares the concern for safe bounds of
Held, Cook and Sewell~\cite{heldcook2012} and of Baldacci et
al.~\cite{baldacci2024}; the smoothed dual is the exponential penalty of
linear programming~\cite{cominetti1994} in the smoothing framework of
Nesterov~\cite{nesterov2005} and of Bienstock and
Iyengar~\cite{bienstock2006}; and the exchange rule extends a theorem of
da Silva et al.~\cite{dasilva2026}. What is new is that these strands
follow from one transform and act as exact tools inside one
branch-and-price algorithm.

The remainder of the paper is organized as follows.
Section~\ref{sec:lehmer} introduces the Lehmer transform and derives from
it the dominance theorem, the Lehmer heuristics and the slack identity.
Section~\ref{sec:entropic} develops entropic column generation.
Section~\ref{sec:perfect} treats perfect packings, with exchange
reductions, blocker branching and backdoor peeling.
Section~\ref{sec:bp} assembles the certified branch-and-price, with the
certification rule, the decision stop and type-split pricing, and proves
its soundness.
Section~\ref{sec:results} reports the computational study. We conclude in
Section~\ref{sec:conclusion}; proofs that are long are collected in
Appendix~\ref{app:proofs}.

\section{Lehmer Transform}\label{sec:lehmer}

This section replaces a bin, a finite multiset of item weights, by an
analytic curve and reads off that curve the dominance order, the
constructive heuristics and a dual solution of the master. With that dual
solution, the accounting identity of column generation becomes an
identity in integers, which tells when the relaxation can decide an
instance and when only a perfect packing can improve the incumbent.

An instance consists of $n$ item types with distinct integer weights $0<w_i\le c$
and demands $d_i\ge1$, to be packed into bins of integer capacity $c$. We
write $D=\sum_id_i$ for the number of items and $W=\sum_id_iw_i$ for their
total weight. A pattern is a nonzero vector $a\in\Z^n_{\ge0}$ with
$\wv^\top a\le c$ and $a\le\dv$; $\calA$ denotes the set of patterns, $N$
their number and $N_{\max}$ the largest number of items in a pattern. The
load and the slack of a pattern are
\begin{equation*}
  \load(a)=\wv^\top a,\qquad \slack(a)=M_1-\load(a),\qquad
  M_1=\max_{a\in\calA}\load(a)\le c,
\end{equation*}
and $M_1=c$ precisely when a full bin exists. With $A$ the matrix whose
columns are the patterns, the pattern formulation of Gilmore and
Gomory~\cite{gilmore1961}, which we call the master, and its linear
relaxation are
\begin{equation}\label{eq:master}
  z^*=\min\bigl\{\one^\top\lambda:\ A\lambda\ge\dv,\ \lambda\in\Z^{N}_{\ge0}\bigr\},\qquad
  \zlp=\min\bigl\{\one^\top\lambda:\ A\lambda\ge\dv,\ \lambda\ge0\bigr\}.
\end{equation}
The covering rows lose nothing, since deleting surplus copies from a
covering solution gives a packing with no more bins. For a dual
vector $\pi\ge0$ the reduced cost of a pattern
and the pricing value are
\begin{equation*}
  \rc_\pi(a)=1-\pi^\top a,\qquad \nu(\pi)=\max_{a\in\calA}\pi^\top a-1,
\end{equation*}
$\pi$ is dual feasible when $\nu(\pi)\le0$, and a feasible $\lambda$ has
over-coverage $e(\lambda)=A\lambda-\dv\ge0$. A restricted master is the
relaxation over a subset of the patterns, with value $\zrmp\ge\zlp$. The
classical lower bounds are $L_1=\lceil W/c\rceil$ and the bound $L_2$ of
Martello and Toth~\cite{martello1990}. An instance has the integer
round-up property (IRUP) if $z^*=\lceil\zlp\rceil$. At a node of a
branching tree with incumbent $\UB$ and some bins already fixed, the
improving target is $T=\UB-1-f$, where $f$ is the number of fixed bins;
it is the number of further bins with which a solution improves the
incumbent.

For a positive signal $\mathbf x=(x_1,\dots,x_k)$ the Mellin transform of
its empirical measure is $\frac1k\sum_ix_i^{\,s-1}$, and the Lehmer
transform is its ratio at consecutive
orders~\cite{ataei2022,ataei2026,ataei_vinogradov2026},
\begin{equation*}
  \Leh_{\mathbf x}(s)=\frac{\sum_{i=1}^{k}x_i^{\,s}}{\sum_{i=1}^{k}x_i^{\,s-1}}
  =\frac{Z_{\mathbf x}(s)}{Z_{\mathbf x}(s-1)}
  =\exp\bigl(K_{\mathbf x}(s)-K_{\mathbf x}(s-1)\bigr),\qquad s\in\R,
\end{equation*}
where $Z_{\mathbf x}(s)=\sum_ix_i^s$ is the partition function of the
signal and $K_{\mathbf x}=\log Z_{\mathbf x}$ its free energy. The escorts
$p_i(u)=x_i^u/Z_{\mathbf x}(u)$ form an exponential family with natural
parameter $u$ and statistic $\log X$, so $K'_{\mathbf x}(u)$ and
$K''_{\mathbf x}(u)$ are the escort mean and variance of $\log X$. Hence
$\Leh_{\mathbf x}$ is real-analytic and, unless the signal is constant,
strictly increasing, the classical monotonicity of the Lehmer
means~\cite{lehmer1971,bullen2003}. The curve runs from $\min_ix_i$ as
$s\to-\infty$ to $\max_ix_i$ as $s\to\infty$, through the harmonic mean at
$s=0$, the arithmetic mean at $s=1$ and the contraharmonic mean at $s=2$;
the order $s$ is the suddency moment of~\cite{ataei2022}, which shifts the
weight of the escort toward the largest observations as it grows. A
pattern $a$ is the signal of its item weights with multiplicities,
\begin{equation*}
  Z_a(s)=\sum_{i=1}^na_iw_i^{s},\qquad K_a(s)=\log Z_a(s),
\end{equation*}
with $Z_a(0)$ its number of items and $Z_a(1)=\load(a)$; by the monotonicity
and the limits of the transform, the order is a dial between the smallest
and the largest items of a bin. The results below use the transform
through its numerator $Z_a$, the partition function of the bin, and
through its free energy $K_a$, whose difference $K_a(s)-K_a(s-1)$ is the
logarithm of the transform; the same construction on the patterns gives
entropic column generation. All Lehmer quantities are
computed in log-space.

The first object read off the curve is the dominance order. Martello and
Toth~\cite{martello1990} call a pattern $F_1$ dominant over $F_2$, both
seen as multisets of item copies, if the items of $F_2$ can be partitioned into blocks $P_1,\dots,P_l$ and $F_1$
contains $l$ distinct items $i_1,\dots,i_l$ with
$w_{i_h}\ge\sum_{k\in P_h}w_k$ for every $h$. Deciding dominance is a
variable-size packing problem; the curves of the two patterns decide it in
one direction by a comparison. Let $\calC$ be the cone of
nondecreasing superadditive functions $f:[0,c]\to\R_{\ge0}$, which vanish
at $0$ since $f(0)\ge2f(0)$; the cone contains
$x\mapsto x^s$ for $s\ge1$ and $x\mapsto\lfloor x/t\rfloor$ for $t>0$; its
elements with $f(c)\le1$ are dual-feasible functions, and every maximal
dual-feasible function~\cite{clautiaux2010} lies in it. To state what the
cone decides, let $v_1,\dots,v_q$ be the distinct weights of $F_2$ and
$m\in\Z^q_{\ge0}$ their multiplicities. A fractional packing of $F_2$ into
the bins $F_1$ is a family of numbers $\lambda_{i,p}\ge0$, for $i\in F_1$
and $p\in\Z^q_{\ge0}$ with $v^\top p\le w_i$, such that
\begin{equation}\label{eq:fracpack}
  \sum_p\lambda_{i,p}\le1\quad\text{for every }i\in F_1,\qquad
  \sum_{i,p}\lambda_{i,p}\,p\ge m ;
\end{equation}
dominance is the same system with $\lambda_{i,p}\in\{0,1\}$.

\begin{theorem}[Dominance]\label{thm:dominance}
  Let $F_1$ and $F_2$ be patterns. Then the following hold.
  \begin{thmlist}
  \item If $F_1$ dominates $F_2$, then
    \begin{equation}\label{eq:dominance}
      \sum_{i\in F_1}f(w_i)\;\ge\;\sum_{k\in F_2}f(w_k)\qquad\text{for every }f\in\calC,
    \end{equation}
    and in particular $K_{F_1}(s)\ge K_{F_2}(s)$ for every $s\ge1$.
  \item Inequality~\eqref{eq:dominance} holds for every $f\in\calC$ if
    and only if $F_2$ has a fractional packing into the bins $F_1$.
  \end{thmlist}
\end{theorem}

\begin{proof}
  \pitem{i} For each block $P_h$, monotonicity and then superadditivity of
  $f$ give
  \begin{equation*}
    f(w_{i_h})\;\ge\;f\Bigl(\sum_{k\in P_h}w_k\Bigr)\;\ge\;\sum_{k\in P_h}f(w_k).
  \end{equation*}
  Summing over $h$ and adding the nonnegative contribution of the items of
  $F_1$ not among $i_1,\dots,i_l$ gives inequality~\eqref{eq:dominance},
  and the power $x\mapsto x^s$ lies in $\calC$ for $s\ge1$.
  \pitem{ii} Suppose first that $\lambda$ satisfies
  system~\eqref{eq:fracpack}, and let $f\in\calC$. Since $f\ge0$, the second
  condition, superadditivity, monotonicity and the first condition give
  \begin{equation*}
    \begin{aligned}
    \sum_{k\in F_2}f(w_k)=\sum_{j}m_jf(v_j)
    &\le\sum_{i,p}\lambda_{i,p}\sum_jp_jf(v_j)
    \le\sum_{i,p}\lambda_{i,p}f(v^\top p)\\
    &\le\sum_{i,p}\lambda_{i,p}f(w_i)
    \le\sum_{i\in F_1}f(w_i).
    \end{aligned}
  \end{equation*}
  Conversely, suppose that system~\eqref{eq:fracpack} has no solution. By
  Farkas' lemma there are $\mu\ge0$, indexed by $F_1$, and $\pi\ge0$, indexed
  by the weights, with $\pi^\top p\le\mu_i$ for every admissible pair $(i,p)$
  and $\pi^\top m>\sum_i\mu_i$. The function
  $f(x)=\max\{\pi^\top p:\ p\in\Z^q_{\ge0},\ v^\top p\le x\}$ on $[0,c]$ is
  nonnegative and nondecreasing, and it is superadditive because $p+p'$ is
  admissible at $x+y$ whenever $p$ is at $x$ and $p'$ at $y$, so $f\in\calC$.
  Moreover $f(w_i)\le\mu_i$ for $i\in F_1$ and $f(v_j)\ge\pi_j$, so
  $\sum_{k\in F_2}f(w_k)\ge\pi^\top m>\sum_i\mu_i\ge\sum_{i\in F_1}f(w_i)$, and
  inequality~\eqref{eq:dominance} fails for $f$.
\end{proof}

Part (i) makes the map $F\mapsto K_F|_{[1,\infty)}$ a monotone
representation of the dominance order, and the restriction to $s\ge1$ is
sharp, since the power is subadditive for $0<s<1$: the pattern $\{4\}$
dominates $\{2,2\}$, while $4^s<2\cdot2^s$ for $s<1$. Part (ii), the
superadditive duality behind dual-feasible functions applied to a pair of
bins, locates exactly what the curves cannot see.

\begin{remark}
  The power family alone does not decide dominance: with $c=10$ the
  patterns $\{7,2\}$ and $\{4,4\}$ satisfy $Z_{\{7,2\}}(s)\ge Z_{\{4,4\}}(s)$
  for all $s\ge1$ without dominance, while $f=\lfloor x/4\rfloor\in\calC$
  separates them. The full cone decides only the fractional relaxation:
  the fifteen items $21,10^{\times3},9^{\times5},7^{\times3},6^{\times3}$
  fill five bins of size $27$ fractionally but not integrally, so in a bin
  of capacity $135$ the pattern of five items of weight $27$ satisfies
  inequality~\eqref{eq:dominance} against the pattern of these items
  without dominating it.
\end{remark}

The curve also defines constructive heuristics. For $s\ge1$ the heuristic
SS-$s$ repeatedly removes from the residual instance the pattern
maximizing $Z_a(s)$, a bounded knapsack with profits $w_i^s$, until all
demand is met. At $s=1$ this is the minimum bin slack heuristic of Gupta
and Ho~\cite{gupta1999} in its exact subset-sum form, which fills each bin
as much as possible; at large orders the largest items dominate the power
sum, and the family reaches, in its limit SS-$\infty$ as $s\to\infty$, the
oldest heuristic of the field.

\begin{theorem}[SS-$\infty$ is first fit decreasing]\label{thm:ssffd}
  Let $D$ be the number of items of the residual and $s>(c+1)\ln D$. The
  maximizer of $Z_a(s)$ over the patterns of the residual is the bin built
  by repeatedly adding the largest residual item that fits, and SS-$s$
  produces the bins of first fit decreasing.
\end{theorem}

The proof is given in Appendix~\ref{app:proofs}. The family thus
interpolates between filling a bin and placing the largest items first,
the two principles behind most bin packing heuristics. For $s>1$ the power
is strictly superadditive, so the chain of inequalities that proves part
(i) of Theorem~\ref{thm:dominance} is strict for $f(x)=x^s$ unless
$F_1=F_2$, and every bin of SS-$s$ with $s>1$ is undominated in the
residual from which it was taken. No order is best on every instance, so
\solver{} runs a portfolio of seven heuristics and stops as soon as an
upper bound meets $\max(L_1,L_2)$: first fit decreasing (FFD) and best fit
decreasing (BFD) run first, then SS-1 and SS-1.5; MBS$'$-1.5 and MBS$'$-2, the
MBS$'$ heuristic of Fleszar and Hindi~\cite{fleszar2002}, which fixes the
largest residual item in the bin before completing it, with the completion
chosen by $Z_a(s)$ at $s=1.5$ and $s=2$; subset-sum-$\tilde n$, which first
places each of the $\tilde n$ items heavier than $c/2$ in a bin of its
own, as in the grouping heuristic of Quiroz-Castellanos et
al.~\cite{quiroz2015}, and then fills these bins and new ones by subset
sums~\cite{caprara2004,caprara2005}. All bins of all runs seed the column
pool.

Column generation needs a dual vector before its first linear program,
and the transform supplies one at every order.

\begin{theorem}[Lehmer duals]\label{thm:lehmerdual}
  For real $s$ let $M_s=\max_{a\in\calA}Z_a(s)$ and $\pi^{(s)}_i=w_i^{s}/M_s$.
  Then the following hold.
  \begin{thmlist}
  \item $\pi^{(s)}$ is dual feasible for the relaxation in
    equation~\eqref{eq:master}, hence $\dv^\top\pi^{(s)}\le\zlp$; for
    $s=1$ this reads $z_0\eqdef W/M_1\le\zlp$, which improves on $W/c$
    whenever no full bin exists.
  \item For $s\ge1$, $\pi^{(s)}$ is the restriction to the item weights of
    $f(x)=x^{s}/M_s\in\calC$, and therefore respects dominance.
  \end{thmlist}
\end{theorem}

\begin{proof}
  \pitem{i} $\pi^{(s)\top}a=Z_a(s)/M_s\le1$ for every $a\in\calA$, and weak
  duality gives $\dv^\top\pi^{(s)}\le\zlp$; $M_1<c$ exactly when no full bin
  exists, in which case $W/M_1>W/c$.
  \pitem{ii} For $s\ge1$, $f(x)=x^s/M_s$ lies in $\calC$, and
  Theorem~\ref{thm:dominance} applies.
\end{proof}

The value of the Lehmer dual as a bound is modest; its role is
structural. It centres the stabilization of the column generation, and it
turns the following identity into an identity in integers.

\begin{theorem}[Budget identity]\label{thm:budget}
  For every $\pi\in\R^n_{\ge0}$ and every $\lambda\ge0$ with $A\lambda\ge\dv$,
  \begin{equation}\label{eq:budget}
    \one^\top\lambda\;=\;\dv^\top\pi\;+\;\sum_{a\in\calA}\lambda_a\,\rc_\pi(a)\;+\;\pi^\top e(\lambda),
  \end{equation}
  and if $\pi$ is dual feasible every term of the two sums is nonnegative.
\end{theorem}

\begin{proof}
  By definition of the over-coverage,
  $\sum_a\lambda_a\pi^\top a=\pi^\top A\lambda=\pi^\top(\dv+e)$, and
  substituting into $\one^\top\lambda=\sum_a\lambda_a\bigl(\rc_\pi(a)+\pi^\top a\bigr)$
  gives identity~\eqref{eq:budget}. If $\pi$ is dual feasible then
  $\rc_\pi\ge0$, and $\pi^\top e\ge0$ since $\pi,e\ge0$.
\end{proof}

The number of bins of any solution is thus the dual objective plus two
nonnegative budgets, one spent on reduced costs and one on over-coverage,
which gives the classical bounds and the fixing rules in one line each.

\begin{corollary}[Bounds and budgets]\label{cor:bounds}
  Let $\pi\ge0$ and $\bar\nu\ge\nu(\pi)$. Then the following hold.
  \begin{thmlist}
  \item $\pi/(1+[\bar\nu]_+)$ is dual feasible, so
    $\zlp\ge\dv^\top\pi/(1+[\bar\nu]_+)$; for the duals of an optimal
    restricted master this is Farley's bound
    $\zlp\ge\zrmp/(1+[\bar\nu]_+)$~\cite{farley1990}.
  \item For every $U\ge\zlp$, $\zlp\ge\dv^\top\pi-U[\bar\nu]_+$, and $\zlp$
    lies in $[\LB_{\rm Lagr},\zrmp]$, where $\LB_{\rm Lagr}$ is the best
    bound of parts (i) and (ii) over the exact pricing calls made so far.
  \item If $\pi$ is dual feasible and $T$ is an integer, every integer
    solution of value at most $T$ satisfies
    $\sum_a\lambda_a\rc_\pi(a)+\pi^\top e\le\beta_T(\pi)\eqdef T-\dv^\top\pi$;
    in particular a pattern with $\rc_\pi(a)>\beta_T(\pi)$ belongs to no
    such solution.
  \end{thmlist}
\end{corollary}

\begin{proof}
  \pitem{i} For every $a$, $\pi^\top a\le1+\nu(\pi)\le1+[\bar\nu]_+$, so the
  scaled vector has every pattern value at most $1$ and weak duality gives
  the bound; for the duals of an optimal restricted master
  $\dv^\top\pi=\zrmp$.
  \pitem{ii} Let $\lambda^*$ be an optimal solution of the relaxation. Since
  $\rc_\pi(a)\ge-[\bar\nu]_+$ for every $a$ and $\pi^\top e\ge0$,
  identity~\eqref{eq:budget} gives
  $\zlp\ge\dv^\top\pi-[\bar\nu]_+\one^\top\lambda^*\ge\dv^\top\pi-U[\bar\nu]_+$. The
  restricted master is a restriction of the master, so $\zrmp\ge\zlp$.
  \pitem{iii} By Theorem~\ref{thm:budget}, an integer solution with
  $\one^\top\lambda\le T$ has
  $\sum_a\lambda_a\rc_\pi(a)+\pi^\top e=\one^\top\lambda-\dv^\top\pi\le\beta_T(\pi)$,
  with every term nonnegative since $\pi$ is dual feasible; a pattern with
  $\lambda_a\ge1$ therefore has $\rc_\pi(a)\le\beta_T(\pi)$.
\end{proof}

Part (ii) is the sandwich that ends column generation, whatever the
pricer reports, and part (iii) with dual-feasible duals of a restricted
master is the reduced-cost fixing of Belov and
Scheithauer~\cite{belov2006}; a partial solution whose fixed columns
exceed the budget can be discarded. With the Lehmer dual,
identity~\eqref{eq:budget} multiplied by $M_1$ counts the slack of the
bins used and the weight of the over-covered items.

\begin{theorem}[Slack identity]\label{thm:slack}
  With $\pi=w/M_1$, every feasible solution $\lambda$ of the relaxation
  satisfies
  \begin{equation}\label{eq:slack}
    \sum_{a\in\calA}\lambda_a\,\slack(a)+\sum_{i=1}^nw_ie_i\;=\;M_1\,\bigl(\one^\top\lambda-z_0\bigr),
  \end{equation}
  an identity in integers when $\lambda$ is integer. Let
  $\sigma_T=M_1(T-z_0)$ for an integer $T$, and $\sigma=\sigma_{\UB-1}$.
  Then the following hold.
  \begin{thmlist}
  \item If $\sigma<0$ then $\UB$ is optimal.
  \item Every improving integer solution uses only patterns with
    $\slack(a)\le\sigma$.
  \item If $\sigma=0$ then every improving solution, integer or fractional,
    uses only patterns of load $M_1$ and has no over-coverage.
  \item $\zlp=z_0$ if and only if the master restricted to the patterns of
    load $M_1$ with equality rows is feasible, and then every optimal
    fractional solution uses only such patterns and has no over-coverage.
  \end{thmlist}
\end{theorem}

\begin{proof}
  With $\pi=w/M_1$ we have $\rc_\pi(a)=\slack(a)/M_1$ and
  $\dv^\top\pi=W/M_1=z_0$; identity~\eqref{eq:budget} multiplied by $M_1$ is
  identity~\eqref{eq:slack}. When $\lambda$ is integer, so is $e$, and every
  term of identity~\eqref{eq:slack} is an integer. All terms of its left
  side are nonnegative for every feasible $\lambda$, and its right side is
  at most $\sigma$ when $\one^\top\lambda\le\UB-1$.
  \pitem{i} A feasible solution with $\one^\top\lambda\le\UB-1$ would make
  the right side negative, so $\zlp>\UB-1$ and every integer solution has
  at least $\UB$ bins.
  \pitem{ii} Every pattern $a$ used by an improving integer solution has
  $\lambda_a\ge1$, so $\slack(a)\le\lambda_a\slack(a)\le\sigma$.
  \pitem{iii} At $\sigma=0$ the right side of an improving solution is at
  most $0$, so every term vanishes: every pattern with $\lambda_a>0$ has
  slack zero, that is, load $M_1$, and $w_i>0$ forces $e=0$.
  \pitem{iv} Applying the identity to an optimal fractional solution, with
  $\one^\top\lambda=\zlp=z_0$, every term vanishes, so the solution is
  supported on patterns of load $M_1$ with $e=0$; conversely a feasible
  solution of the equality master over the patterns of load $M_1$ has value
  $W/M_1=z_0$, so $\zlp\le z_0\le\zlp$.
\end{proof}

At a node the theorem applies to the residual instance, with the
improving target $T$ in place of $\UB-1$, and $\sigma$ then denotes
$\sigma_T$. The integer $\sigma$ is the total slack that an improving
packing may spend: it bounds the slack of the columns fixed in a partial
solution and restricts the pricing of a node to the loads in
$[M_1-\sigma,M_1]$. It also separates two regimes. When $\sigma>0$ the
relaxation guides the search. When $\sigma=0$, an improving solution
fills every bin to the load $M_1$, so it is a perfect packing for the
capacity $M_1$, and a perfect packing in the usual sense when $M_1=c$. On
the AI and ANI classes this holds with $M_1=c$ at
every root once first fit decreasing returns $L_1+1$ bins, since $W=L_1c$
and full bins exist by construction; as moreover $\zlp=z_0=L_1$ there,
the relaxation is a fractional perfect packing by part (iv), every full
pattern is tight at the Lehmer dual, and the relaxation cannot tell which
full patterns form an integer perfect packing.

\section{Entropic Column Generation}\label{sec:entropic}

This section applies the construction of the Lehmer transform to the
patterns. Weighted by the exponential of its value under a dual vector,
each pattern becomes a state of a probability distribution whose
partition function the knapsack recursion computes once maximization is
replaced by summation. Its free energy brackets the pricing value and
smooths the Lagrangian dual of the master, which gives bounds at every
iteration of column generation, and sampling from it gives restricted
masters of bounded size. In \solver{} the ensemble is a source of bounds,
a measure of degeneracy and a seed for the root column generation, while
exact pricing drives the column generation itself.

Fix $\pi\ge0$, write $v_a=\pi^\top a$ for the value of a pattern and
$v^*=\max_av_a=1+\nu(\pi)$ for the largest value. For an order $s\ge0$ the
free energy of the pattern set, the Gibbs distribution on the patterns and
its entropy are
\begin{equation*}
  F(s)=\log\sum_{a\in\calA}e^{sv_a},\qquad
  P_s(a)=e^{sv_a-F(s)},\qquad
  H(P_s)=-\sum_{a\in\calA}P_s(a)\log P_s(a).
\end{equation*}
With the values $x_a=e^{v_a}$, $F$ is the free energy of these values and
$F(s)-F(s-1)$ the logarithm of their Lehmer transform: the ensemble is
the Lehmer transform at the scale of patterns. Its partition function is
the bounded knapsack recursion over loads with maximization replaced by
summation, in log-space with $u\oplus v=\log(e^u+e^v)$,
\begin{equation*}
  G_i[C]=\bigoplus_{k=0}^{\min(d_i,\lfloor C/w_i\rfloor)}\bigl(G_{i-1}[C-kw_i]+sk\pi_i\bigr),
  \qquad F(s)=\bigoplus_{C=1}^{c}G_n[C],
\end{equation*}
which we call the sum--product knapsack; the load $C=0$, which holds only
the zero vector, is left out, and binary splitting of the demands, under
which count vectors are not represented uniquely, is not admissible. A
pass costs
$O(c\sum_i\min(d_i,c/w_i))$, the order of one exact pricing pass; a
backward pass gives the marginals $\E_{P_s}[a_i]$, forward filtering with
backward sampling draws exact samples in
$O\bigl(\log c+\sum_i\min(d_i,\lfloor c/w_i\rfloor)\bigr)$ each after an $O(c)$
preparation, and the pass at
$s=0$ counts the patterns, $N=e^{F(0)}$. As for every exponential family,
$F$ is convex with $F'(s)=\E_{P_s}[v]$ and $F''(s)=\Var_{P_s}(v)$, and
\begin{equation}\label{eq:gibbsid}
  F(s)=s\,\E_{P_s}[v]+H(P_s),
\end{equation}
with $H(P_0)=\log N$ and $H(P_s)\to\log N^*$ as $s\to\infty$, where $N^*$
counts the patterns of value $v^*$. Identity~\eqref{eq:gibbsid} is the
source of every bound below.

\begin{theorem}[Entropic bounds]\label{thm:entbounds}
  For every $s>0$ the following hold.
  \begin{thmlist}
  \item $\E_{P_s}[v]\le v^*\le F(s)/s$, and both gaps are at most
    $H(P_s)/s\le(\log N)/s$.
  \item $F(s)/s$ is nonincreasing and $\E_{P_s}[v]$ nondecreasing in $s$,
    and both tend to $v^*$.
  \item If $\Delta$ is the gap between $v^*$ and the next largest pattern
    value, with $\Delta=+\infty$ when $N^*=N$, then
    \begin{equation}\label{eq:gapbound}
      \frac{F(s)}{s}-v^*\;\le\;\frac1s\log\bigl(N^*+(N-N^*)e^{-s\Delta}\bigr).
    \end{equation}
  \item With $\bar\nu_s=F(s)/s-1$ and every $U\ge\zlp$,
    \begin{equation*}
      \zlp\;\ge\;\frac{\dv^\top\pi}{1+[\bar\nu_s]_+}\qquad\text{and}\qquad
      \zlp\;\ge\;\dv^\top\pi-U[\bar\nu_s]_+ .
    \end{equation*}
  \end{thmlist}
\end{theorem}

\begin{proof}
  The family $\{P_s\}$ is the exponential family with natural parameter
  $s$, statistic $v$ and cumulant $F$, which gives $F'=\E_{P_s}[v]$ and
  $F''=\Var_{P_s}(v)\ge0$; from $\log P_s(a)=sv_a-F(s)$,
  $H(P_s)=-\sum_aP_s(a)\bigl(sv_a-F(s)\bigr)=F(s)-sF'(s)$, which is
  identity~\eqref{eq:gibbsid}. At $s=0$ the distribution is uniform on
  $\calA$, so $H(P_0)=\log N$; as $s\to\infty$ it concentrates uniformly on
  the $N^*$ maximizers.
  \pitem{i} From $e^{sv^*}\le\sum_ae^{sv_a}=e^{F(s)}$ we get $v^*\le F(s)/s$.
  Identity~\eqref{eq:gibbsid} with $\E_{P_s}[v]\le v^*$ gives
  \begin{equation*}
    \begin{aligned}
    \frac{F(s)}{s}-v^*&=\E_{P_s}[v]-v^*+\frac{H(P_s)}{s}\;\le\;\frac{H(P_s)}{s},\\
    v^*-\E_{P_s}[v]&=\frac{sv^*-F(s)}{s}+\frac{H(P_s)}{s}\;\le\;\frac{H(P_s)}{s},
    \end{aligned}
  \end{equation*}
  the second using $F(s)\ge sv^*$, and $H(P_s)\le\log N$.
  \pitem{ii} By identity~\eqref{eq:gibbsid},
  $\frac{d}{ds}\frac{F(s)}{s}=\frac{sF'(s)-F(s)}{s^2}=-\frac{H(P_s)}{s^2}\le0$,
  and $F'=\E_{P_s}[v]$ is nondecreasing because $F''\ge0$; both tend to $v^*$
  since the weight of the non-maximizers under $P_s$ is at most
  $(N-N^*)e^{-s\Delta}/N^*$.
  \pitem{iii} Splitting the sum into the maximizers and the rest gives
  $\sum_ae^{sv_a}\le N^*e^{sv^*}+(N-N^*)e^{s(v^*-\Delta)}$; take logarithms and
  divide by $s$.
  \pitem{iv} By part (i), $\nu(\pi)=v^*-1\le\bar\nu_s$, and parts (i) and (ii)
  of Corollary~\ref{cor:bounds} apply with $\bar\nu=\bar\nu_s$.
\end{proof}

One pass of the sum--product knapsack therefore brackets the pricing
value, and by inequality~\eqref{eq:gapbound} the upper bracket is within
about $(\log N^*)/s$ of $v^*$ once
$s\Delta\gtrsim\ln\bigl((N-N^*)/N^*\bigr)$. Part (iv) gives a lower bound
on the relaxation at every iteration, from any dual vector and without
exact pricing, while the restricted master bounds it from above. At an
optimal restricted master the basic columns have value one, so
$F(s)/s>1$ and no entropic quantity can certify optimality; that remains
the task of exact pricing, and the pass runs in double precision.

The Lagrangian dual of the master can be written with the pricing value
as a penalty. For $U\ge\zlp$ let
\begin{equation*}
  \theta(\pi)=\dv^\top\pi-U[\nu(\pi)]_+ .
\end{equation*}
The penalty is not smooth, since $\nu$ is a maximum of linear functions,
and the free energy of the ensemble on $\calA\cup\{0\}$ smooths it in the
sense of Nesterov~\cite{nesterov2005}, with a loss that the entropy
controls; this is the exponential penalty of linear programming studied by
Cominetti and San Mart\'in~\cite{cominetti1994}, applied to the pattern
master. For $s>0$ let
\begin{equation*}
  g_s(\pi)=\frac1s\log\Bigl(1+\sum_{a\in\calA}e^{s(v_a-1)}\Bigr),\qquad
  \theta_s(\pi)=\dv^\top\pi-Ug_s(\pi),
\end{equation*}
and let $\tilde P_s$ be the distribution on $\calA\cup\{0\}$ with
$\tilde P_s(a)\propto e^{s(v_a-1)}$ for $a\in\calA$ and $\tilde P_s(0)\propto1$.

\begin{theorem}[Smoothed dual]\label{thm:smoothed}
  Let $U\ge\zlp$ and $s>0$. Then the following hold.
  \begin{thmlist}
  \item $\theta\le\zlp$ on $\R^n_{\ge0}$, and $\max_{[0,1]^n}\theta=\zlp$,
    attained at every optimal dual; if $U>\zlp$, $\theta$ is coercive on
    $\R^n_{\ge0}$ and its maximizers there are exactly the optimal duals.
  \item $\theta-U\log(1+N)/s\le\theta_s\le\theta$.
  \item $\theta_s$ is concave and $C^\infty$, with gradient
    $\dv-U\E_{\tilde P_s}[a]$ that is Lipschitz with constant
    $L_s=UsN_{\max}^2$.
  \item Nesterov's accelerated projected gradient method on $-\theta_s$
    over $[0,1]^n$ with step $1/L_s$ satisfies
    \begin{equation}\label{eq:nesterov}
      \max_{[0,1]^n}\theta_s-\theta_s(\pi^k)\;\le\;\frac{2UsN_{\max}^2n}{(k+1)^2},
    \end{equation}
    every iterate carrying the valid bound $\zlp\ge\theta_s(\pi^k)$, and
    with
    \begin{equation}\label{eq:nesterovparams}
      s=\frac{2U\log(1+N)}{\epsilon}\qquad\text{and}\qquad
      k+1\ge\frac{2\sqrt2\,UN_{\max}\sqrt{n\log(1+N)}}{\epsilon}
    \end{equation}
    one has $\zlp-\epsilon\le\theta_s(\pi^k)\le\zlp$.
  \end{thmlist}
\end{theorem}

The proof is given in Appendix~\ref{app:proofs}. Part (iv) gives an
$\epsilon$-certificate of $\zlp$, as a worst-case guarantee, after
$O(UN_{\max}\sqrt{n\log N}/\epsilon)$ passes of the sum--product knapsack
and no linear program, a rate linear in $1/\epsilon$ against
$1/\epsilon^2$ for multiplicative-weights schemes with the same
oracle~\cite{garg2007}, as for the smoothing method of Bienstock and
Iyengar~\cite{bienstock2006}. What is specific to cutting stock is that
the smoothed oracle costs the order of one exact pricing pass, and that
the maximizer of the smoothed dual carries a primal solution, its ensemble
scaled by the cap $U$, with a duality gap that is an entropy.

\begin{theorem}[Entropic duality]\label{thm:entduality}
  Let $U>\zlp$ and $s>0$, let $\hat\pi$ maximize $\theta_s$ over
  $\R^n_{\ge0}$, which exists, let $\tilde P=\tilde P_s(\hat\pi)$, and set
  $\hat\lambda_a=U\tilde P(a)$ for $a\in\calA$. Then the following hold.
  \begin{thmlist}
  \item $\hat\lambda$ is feasible for the relaxation, with equality in
    every row $i$ with $\hat\pi_i>0$.
  \item The duality gap of the pair is the entropy per unit of order,
    \begin{equation}\label{eq:entgap}
      \one^\top\hat\lambda-\theta_s(\hat\pi)=\frac{U\,H(\tilde P)}{s}.
    \end{equation}
  \item $\theta_s(\hat\pi)\le\zlp\le\one^\top\hat\lambda$, and both are
    within $U\log(1+N)/s$ of $\zlp$.
  \end{thmlist}
\end{theorem}

\begin{proof}
  A maximizer exists because $\theta_s\le\theta$ and $\theta$ is coercive
  on $\R^n_{\ge0}$ by part (i) of Theorem~\ref{thm:smoothed}, so the
  superlevel sets of the continuous $\theta_s$ are compact.
  \pitem{i} The first-order conditions of the concave $\theta_s$ over the
  orthant read
  \begin{equation}\label{eq:foc}
    \partial_i\theta_s(\hat\pi)=d_i-U\,\E_{\tilde P}[a_i]\;\le\;0,
    \qquad\text{with equality whenever }\hat\pi_i>0,
  \end{equation}
  that is, $A\hat\lambda\ge\dv$ with equality in the rows with $\hat\pi_i>0$.
  \pitem{ii} Multiplying the equalities in condition~\eqref{eq:foc} by
  $\hat\pi_i$ and summing gives $\dv^\top\hat\pi=U\,\E_{\tilde P}[\hat\pi^\top a]$,
  where the pattern $0$ contributes $0$. Identity~\eqref{eq:gibbsid} for the
  ensemble on $\calA\cup\{0\}$ with the statistic $u_a=v_a-1$ for $a\in\calA$
  and $u_0=0$ reads $sg_s(\hat\pi)=s\,\E_{\tilde P}[u]+H(\tilde P)$. Hence
  \begin{equation*}
    \theta_s(\hat\pi)=\dv^\top\hat\pi-Ug_s(\hat\pi)
    =U\,\E_{\tilde P}\bigl[\hat\pi^\top a-u_a\bigr]-\frac Us H(\tilde P),
  \end{equation*}
  and $\hat\pi^\top a-u_a$ equals $1$ for $a\in\calA$ and $0$ for $a=0$, so the
  expectation is $\tilde P(\calA)=\one^\top\hat\lambda/U$, which gives
  identity~\eqref{eq:entgap}.
  \pitem{iii} $\theta_s(\hat\pi)\le\theta(\hat\pi)\le\zlp$ by
  Theorem~\ref{thm:smoothed}, $\zlp\le\one^\top\hat\lambda$ by part (i), and
  the two differ by $UH(\tilde P)/s\le U\log(1+N)/s$ by part (ii).
\end{proof}

The entropy of the ensemble is thus exactly the price of smoothing. The
primal solution is never written down, since its marginals are what the
backward pass returns, and a restricted master of bounded size is
obtained by sampling.

\begin{theorem}[Sampled master]\label{thm:sampled}
  Let $\tilde P=\tilde P_s(\pi)$ satisfy $U\E_{\tilde P}[a_i]\ge d_i$ for all
  $i$, draw $a^{(1)},\dots,a^{(m)}$ independently from $\tilde P$, and set
  $\hat\lambda^m_a=(U/m)\,\#\{j:a^{(j)}=a\}$. For $\xi,\zeta\in(0,1)$,
  $d_{\min}=\min_id_i$ and
  \begin{equation}\label{eq:msamples}
    m\;\ge\;\frac{2N_{\max}U\ln(2n/\zeta)}{\xi^2d_{\min}},
  \end{equation}
  with probability at least $1-\zeta$ the vector
  $\lambda'=\hat\lambda^m/(1-\xi)$ is feasible, has at most $m$ nonzero
  entries, and satisfies
  \begin{equation}\label{eq:sampledvalue}
    \one^\top\lambda'\;\le\;\frac{U\tilde P(\calA)+U\sqrt{\ln(2/\zeta)/2m}}{1-\xi}.
  \end{equation}
\end{theorem}

The proof is given in Appendix~\ref{app:proofs}. With $\tilde P$ from
Theorem~\ref{thm:entduality}, the restricted master over the samples has,
with probability at least $1-\zeta$, value at most
$\bigl(\zlp+U\log(1+N)/s+U\sqrt{\ln(2/\zeta)/2m}\bigr)/(1-\xi)$ with a
number of columns independent of $N$. Since a pass costs time
proportional to the capacity, coarsening helps. With an integer $1\le g\le\min_iw_i$, every pattern of the
instance is a pattern of the instance with weights $\lfloor w_i/g\rfloor$
and capacity $\lfloor c/g\rfloor$, since
$\sum_ia_i\lfloor w_i/g\rfloor\le\lfloor\load(a)/g\rfloor$, so the upper
brackets computed there remain valid at a $g$-th of the cost.

The ensemble also measures what the relaxation leaves undecided. At an
optimal dual, $N^*$ counts the columns among which the relaxation is
indifferent, and at a finite order $e^{H(P_s)}$ is the effective number of
nearly optimal patterns. Where $\zlp=z_0$ the Lehmer dual is optimal and
its tight patterns are all those of load $M_1$, so the relaxation is
indifferent among all of them; it is at such nodes that combinatorial
reasoning must take over.

\section{Perfect Packings}\label{sec:perfect}

When $\sigma=0$, an improving solution is a perfect packing for the
capacity $M_1$, and the relaxation carries no information about which
patterns of load $M_1$ form one. This section decides such instances with combinatorial tools: reductions fix bins that some perfect packing must
contain, a complete search branches where they stall, and backdoor
peeling holds back a few items so that the reductions can cascade. They
act on the instance at the root of the branch-and-price and on the
residuals of their own search.

Throughout, $c$ stands for the load $M_1$, which equals the capacity on
the instances of interest; a full bin is a bin of load $c$, $R$ is the
multiset of items of the instance or of a residual, a perfect packing of
$R$ is a partition of $R$ into full bins, and a
completion of an item $i\in R$ is a sub-multiset of $R\setminus\{i\}$ of
weight $c-w_i$. An item is large if $2w_i>c$; two large items never share
a bin, and a completion of a large item contains no large item. A bin can
be fixed whenever every perfect packing can be rearranged to contain it,
and an exchange argument tests this with one subset-sum computation.

\begin{lemma}[Exchange]\label{lem:exchange}
  Let $P=\{i\}\cup Q$ be a full bin of items of $R$, and suppose that every
  full bin $F\subseteq R$ with $i\in F$ satisfies $|Q\setminus F|\le1$. If $R$
  has a perfect packing, then it has one that contains the bin $P$.
\end{lemma}

\begin{proof}
  We show more: for a perfect packing whose bin $B$ of $i$ satisfies
  $|Q\setminus B|\le1$, an exchange between $B$ and one other bin yields a
  perfect packing that contains $P$. If $Q\subseteq B$, then $P\subseteq B$ and
  both are full, so $B=P$. Otherwise $Q\setminus B=\{x\}$, the item $x$ lies
  in a bin $B'\ne B$, and $P\setminus\{x\}\subseteq B$. The items
  $R'=B\setminus(P\setminus\{x\})$ weigh $c-(c-w_x)=w_x$, so exchanging $R'$
  and $x$ turns $B$ into $P$ and $B'$ into $(B'\setminus\{x\})\cup R'$, both
  full, and leaves every other bin unchanged. Under the hypothesis of the
  lemma the bin of $i$ in every perfect packing is a full bin containing
  $i$, so $|Q\setminus B|\le1$ holds.
\end{proof}

For $|Q|=1$ this is the classical rule that pairs an item with an exact
complement. For $|Q|=2$ the hypothesis says that $i$ has no completion
avoiding both items of $Q$, a single subset-sum test; in this form the
lemma is Theorem~3 of da Silva et al.~\cite{dasilva2026}, whose ANI
algorithm applies it to the triplets of the construction. The general form
holds for any $Q$ and is tested pair by pair: for every pair $\{x,y\}\subseteq Q$,
no subset of $R\setminus\{i,x,y\}$ may weigh $c-w_i$, one subset-sum test each;
the solver uses it with two and three items. Three further statements hold in every perfect
packing and therefore fix or forbid bins outright.

\begin{lemma}[Forced bins]\label{lem:forced}
  Let $i\in R$. Then the following hold.
  \begin{thmlist}
  \item If $i$ has no completion, $R$ has no perfect packing.
  \item If every completion of $i$ contains an item of weight $w$, then in
    every perfect packing $i$ shares its bin with an item of weight $w$,
    and the two can be merged into one item of weight $w_i+w$.
  \item If $i$ has exactly one completion, as a multiset of weights, the
    bin it forms with $i$ belongs to every perfect packing.
  \end{thmlist}
\end{lemma}

\begin{proof}
  In a perfect packing the bin $B$ of $i$ is full, so $B\setminus\{i\}$ is a
  completion of $i$.
  \pitem{i} No completion means no bin for $i$.
  \pitem{ii} $B\setminus\{i\}$ contains an item of weight $w$. Items of equal
  weight are interchangeable, so the perfect packings in which $i$ shares a
  bin with an item of weight $w$ are those of the multiset in which the two
  are replaced by one item of weight $w_i+w$.
  \pitem{iii} $B\setminus\{i\}$ is the one completion, up to the interchange of
  items of equal weight.
\end{proof}

Lemma~\ref{lem:exchange} preserves the existence of a perfect packing and
Lemma~\ref{lem:forced} holds in every one, so applying them in any order
until none applies leaves a residual that has a perfect packing if and
only if $R$ has one; a perfect packing of the residual, with the fixed
bins and the merged items expanded, is a perfect packing of $R$. The tests
run on the subset-sum counts of the residual modulo $61$-bit primes.
Removing an item of weight $w$ divides the generating polynomial by
$1+z^w$ in time $O(c)$, and for a full triple $\{i,x,y\}$ the number of
completions of $i$ that avoid $x$ and $y$ is
\begin{equation}\label{eq:avoidcount}
  N_i-N_{i,x}-N_{i,y}+1,
\end{equation}
where $N_i$ counts the completions of $i$ and $N_{i,x}$ those containing a
given copy of $x$; once the counts without $i$ are computed, $N_i$ is one
of them and $N_{i,x}$ takes time proportional to $c/w_x$. The modular
counts serve only as a filter, since a residue says nothing about a count
that may be astronomically large. A zero, on which
Lemma~\ref{lem:exchange} and parts (i) and (ii) of Lemma~\ref{lem:forced}
rest, is confirmed by an exact bitset reachability test before the rule
is applied. A single completion $C$, on which part (iii) rests, is
confirmed by one reachability test for each weight of $C$, showing that no
completion uses fewer copies of that weight; any other completion would,
since a completion that contains $C$ as a multiset equals it. Every
applied rule is logged with its test, so that the chain of reductions is
a certificate that independent code can re-verify.

When no reduction applies, the question is decided by a search whose
nodes are residuals, with the reductions run to a fixpoint at every node;
a node of the search is a plain instance with no branching rule attached,
so Lemmas~\ref{lem:exchange} and~\ref{lem:forced} remain valid there. The
branching rule exploits why an exchange fails: if a full triple fails the
test of Lemma~\ref{lem:exchange} only because of a few items, those items
are the ones to branch on.

\begin{lemma}[Blocker branching]\label{lem:blocker}
  Let $\{i,x,y\}\subseteq R$ be a full bin and $S\subseteq R\setminus\{i,x,y\}$
  a set such that $i$ has no completion in $R\setminus(\{i,x,y\}\cup S)$. Then
  the following hold.
  \begin{thmlist}
  \item If $R$ has a perfect packing in which the bin of $i$ contains no item
    of $S$, it has one that contains the bin $\{i,x,y\}$.
  \item If $R$ has a perfect packing, then so does $R$ with the bin
    $\{i,x,y\}$ fixed, or $R$ with $i$ and some $s\in S$ merged.
  \end{thmlist}
\end{lemma}

\begin{proof}
  \pitem{i} Let $B$ be the bin of $i$, with $B\cap S=\emptyset$. Its completion
  $B\setminus\{i\}$ lies in $R\setminus(\{i\}\cup S)$ and, by hypothesis, not in
  $R\setminus(\{i,x,y\}\cup S)$, so $B$ contains $x$ or $y$. Hence
  $|\{x,y\}\setminus B|\le1$, and the exchange in the proof of
  Lemma~\ref{lem:exchange} gives a perfect packing that contains $\{i,x,y\}$.
  \pitem{ii} In a perfect packing the bin of $i$ either contains no item of
  $S$, and part (i) applies, or contains some $s\in S$, and then $i$ and $s$
  share a bin and can be merged as in part (ii) of Lemma~\ref{lem:forced}.
\end{proof}

The search therefore branches into the child that fixes $\{i,x,y\}$ and
one child per item of the blocker set $S$, choosing the triple with the
smallest $S$ and, among those, the largest item $i$ first; blocker sets of
one or two items are found by testing the items of one or two conflicting
completions. Where no such set exists it branches on the item with the
fewest completions, as ranked by the subset-sum counts, with one child per
completion; the completions are enumerated exactly, as multisets of
weights, since a missing one would make the branching incomplete.
Residuals that failed are remembered by a hash of their weight multiset. A
packing found closes the instance after a check bin by bin; an exhausted
search proves that no perfect packing exists, and its whole tree, the
blocker sets, the enumerated completions and every residual pruned by the
memory included, is written to the certificate. The certificate is
re-verified by a separate checker with its own subset-sum dynamic program
and its own enumeration of completions, which compares every memoised
residual with the stored multiset in full, so that a hash collision cannot
enter a proof.

The cascade of exchange reductions can stall at its first step because a
handful of items give the large items completions they would not
otherwise have. Satisfiability solving calls a set of variables a weak
backdoor if some assignment to them leaves a formula that propagation
solves, and a strong backdoor if every assignment does~\cite{williams2003};
in integer programming a backdoor is a small set of variables on which
branching alone solves the problem~\cite{dilkina2009,fischetti2013}, and
such sets can be found by search~\cite{khalil2022}. The analogue here is
of the weak kind: a small multiset $V\subseteq R$ of items held back from
the bins of the large items, under which hypothesis the completions of a
large item are computed in $R\setminus V$ and the reductions may run
again. It serves to find a packing and never to refute one.

\begin{lemma}[Conditional exchange]\label{lem:condexchange}
  Let $V\subseteq R$, let $J$ be a set of large items of $R$ disjoint from
  $V$, and say that a perfect packing respects $(V,J)$ if the bins of the
  items of $J$ contain no item of $V$. Let $i\in J$. Then the following
  hold.
  \begin{thmlist}
  \item Let $\{i,x,y\}$ be a full bin with $x,y\notin V$ such that every full
    bin $F\ni i$ with $F\cap V=\emptyset$ contains $x$ or $y$. If $R$ has a
    perfect packing that respects $(V,J)$, it has one that respects $(V,J)$
    and contains the bin $\{i,x,y\}$.
  \item If $i$ has exactly one completion in $R\setminus V$, as a multiset of
    weights, that completion forms the bin of $i$ in every perfect packing
    that respects $(V,J)$.
  \item Let $x\notin V$ be an item with $w_i+w_x=c$. If $R$ has a perfect
    packing that respects $(V,J)$, it has one that respects $(V,J)$ and
    contains the bin $\{i,x\}$.
  \end{thmlist}
\end{lemma}

The proof, which follows the exchanges of Lemma~\ref{lem:exchange} and
checks that they keep the bins of $J$ free of $V$, is given in
Appendix~\ref{app:proofs}. We call backdoor peeling the procedure that
learns a reserve $V$ and runs the reserved cascade under it. The cascade
takes for $J$ the large items that have a completion in $R\setminus V$,
applies parts (i) and (iii) to them, larger weights first, and part (ii)
when it stalls. A large item that loses its last completion in
$R\setminus V$ cannot respect the hypothesis, so it leaves $J$ and is
exempt from the conditional rules, and a reserve that exempts too many is
rejected. The cascade stops as soon as the residual has at most $\kappa$
items, so that the few large items that belong to the residual's own bins
are never fixed under a hypothesis they may violate. The residual is then
decided by the search above without any hypothesis, and a packing found,
with the fixed bins, is checked bin by bin. The hypothesis is never proved
and never needs to be: a wrong reserve costs time, and every packing
returned is verified on the original items.

The reserve is learned from the stalls of the cascade. At a stall, a full
triple $\{i,x,y\}$ whose test fails because of few conflicting
completions, between one and $\mu$ of them by count~\eqref{eq:avoidcount}
on $R\setminus V$, names its own obstruction: any reserve that makes each
conflicting completion impossible lets the cascade fix $\{i,x,y\}$. The
candidates are the sets of at most four item copies that do so, ranked by
the number of conflicts, together with the analogous sets for the large
items with the fewest completions; their items belong to completions of
large items, so they are below half the capacity and $V$ never meets $J$.
Each candidate is scored by the number of bins its cascade fixes, and the
best few are also rolled out: their score includes that of their best
continuation one stall further, because the items of a backdoor can open
the cascade only together, one of them showing its effect only after the
other has been held back. The search over reserves is best-first with a
beam of width $b$ and at most $n_V$ held-back items.

\section{Certified Branch-and-Price}\label{sec:bp}

A branch-and-price that reports proved optima must make exact the bounds
by which it prunes and the moment at which the column generation of a
node stops. This section does both, prices the branching rules and cuts
without a state per cut, and assembles the solver.

A bound proves optimality only if it survives the rounding errors of the
computation that produced it. Pricing therefore runs on scaled integer
profits: with a scale $2^k$, $p_i=\lfloor2^k\pi_i\rfloor$ and
$P^*=\max_{a\in\calA}p^\top a$ computed exactly on integers,
$2^k\pi_i\le p_i+1$ along at most $N_{\max}$ items gives
\begin{equation}\label{eq:scaling}
  \nu(\pi)\;\le\;\bar\nu\;\eqdef\;\frac{P^*+N_{\max}}{2^k}-1 .
\end{equation}
A pricing pass that finds no pattern of negative reduced cost therefore
proves $\nu\le\bar\nu\approx N_{\max}2^{-k}$, not $\nu\le0$, and only the
lower end of the sandwich of part (ii) of Corollary~\ref{cor:bounds},
computed with $\bar\nu$, is used as a bound, so that no tolerance of the
linear programming solver enters a proof. The solver may return duals
with slightly negative components. Deleting from a pattern its items of
negative scaled profit leaves a pattern of no smaller value or no item at
all, and every item forms a pattern alone, so the pricing maximum of the
profits $\max(p,0)$ is at most $\max(P^*,0)$, and $P^*<0$ only if
$\pi\le0$. Hence $\nu(\pi^+)\le[\bar\nu]_+$ for $\pi^+=\max(\pi,0)$, so
Corollary~\ref{cor:bounds} applies to $\pi^+$ with $[\bar\nu]_+$ in place of
$\bar\nu$, and since $\dv^\top\pi\le\dv^\top\pi^+$ its bounds computed from $\pi$
remain valid. The scale limits what such a bound can certify.

\begin{lemma}[Resolution of a certified bound]\label{lem:resolution}
  Let $\pi\ge0$, let $\bar\nu$ be given by inequality~\eqref{eq:scaling}
  at scale $2^k$, and let $T\ge0$ and $\tau>0$. Then the following hold.
  \begin{thmlist}
  \item $\nu(\pi)\le\bar\nu\le\nu(\pi)+N_{\max}2^{-k}$, with equality on
    the right when $2^k\pi$ is integral.
  \item The Farley bound $\dv^\top\pi/(1+[\bar\nu]_+)$ exceeds $T+\tau$ if
    and only if
    \begin{equation}\label{eq:precision}
      [\bar\nu]_+\;<\;\frac{\dv^\top\pi-T-\tau}{T+\tau}.
    \end{equation}
  \item If $\nu(\pi)\le\varepsilon$ for some $\varepsilon\ge0$, the Farley
    bound exceeds $T+\tau$ whenever
    $\dv^\top\pi-T>\tau+(T+\tau)(\varepsilon+N_{\max}2^{-k})$; if $2^k\pi$
    is integral and $\nu(\pi)=0$, it does so if and only if
    $\dv^\top\pi-T>\tau+(T+\tau)N_{\max}2^{-k}$.
  \end{thmlist}
\end{lemma}

\begin{proof}
  \pitem{i} The left inequality is inequality~\eqref{eq:scaling}. For the
  right one, $p_i\le2^k\pi_i$ gives
  $P^*\le2^k\max_a\pi^\top a=2^k(1+\nu(\pi))$, hence
  $\bar\nu\le\nu(\pi)+N_{\max}2^{-k}$, with equality when $p=2^k\pi$.
  \pitem{ii} The Farley bound exceeds $T+\tau$ if and only if
  $\dv^\top\pi>(T+\tau)(1+[\bar\nu]_+)$, which is
  condition~\eqref{eq:precision} since $T+\tau>0$.
  \pitem{iii} If $\nu(\pi)\le\varepsilon$, part (i) gives
  $[\bar\nu]_+\le\varepsilon+N_{\max}2^{-k}$, and the stated lift makes the
  right side of condition~\eqref{eq:precision} exceed
  $\varepsilon+N_{\max}2^{-k}$. If $2^k\pi$ is integral and $\nu(\pi)=0$, part
  (i) gives $[\bar\nu]_+=N_{\max}2^{-k}$, and condition~\eqref{eq:precision}
  becomes the stated lift.
\end{proof}

By part (iii), a lift $\dv^\top\pi-T$ of the order of $TN_{\max}2^{-k}$ is
the smallest that scale $2^k$ can certify in general. A coarse scale is
fast but cannot certify the tiny lifts that cuts produce on the hardest
instances, so the solver prices at a coarse scale $k_0$ and recomputes
every bound that must certify a small lift by one exact pass at a fine
scale $k_1$. The pricing maximum is exact, and the bound itself is
evaluated in floating point; by the standard error analysis of
floating-point summation~\cite{higham2002}, its error is at most about
$N\mathrm{u}$ times the sum of the absolute values of its $N$ terms, with
$\mathrm{u}=2^{-53}$ the unit roundoff. The solver subtracts this error
bound from the computed value at every evaluation, so that the bound it
uses never exceeds the exact one. The bound of a node is
\begin{equation}\label{eq:lbnode}
  \LB_{\rm node}=\max\bigl\{L_1,\ L_2,\ f+\lceil\LB_{\rm Lagr}-\tau\rceil\bigr\},
\end{equation}
a bound on the total number of bins of every solution of the node, with
$f$ its number of fixed bins and $\LB_{\rm Lagr}$ computed from $\bar\nu$,
and never $\lceil\zrmp\rceil$, which is one too many when $\zlp$ is an
integer and the restricted master stops slightly above it; at a node,
$\zlp$ and $\zrmp$ denote the values of the relaxation and of the
restricted master over the node's pattern space and rows. A node is pruned
when $\LB_{\rm node}\ge\UB$, and a child inherits the bound of its parent,
which holds for it since its solutions are solutions of the parent. One
rule governs every closure.

\begin{definition}[Certification rule]\label{def:certrule}
  A node, the root included, is closed only in one of the following ways.
  \begin{thmlist}
  \item By a packing of the node, checked bin by bin, whose number of bins
    attains the bound~\eqref{eq:lbnode} of the node.
  \item By a certified bound: $\LB_{\rm node}\ge\UB$, or
    $\lceil\LB-\tau\rceil>T$, where $\LB$ is a bound of part (i) or (ii) of
    Corollary~\ref{cor:bounds}, the latter with a cap $U\ge T+1$ and
    with the duals clipped to their signs, evaluated with a value
    $\bar\nu\ge\nu(\pi)$ returned by an integer pricing pass over the node's
    pattern space.
  \item By an integer test: $\sigma<0$, an item that no pattern of the
    node's space covers, or a Farkas ray $y\ge0$ of an infeasible restricted
    master over all its rows, rounded down to integers, whose certificate an
    exact pricing pass over the node's pattern space confirms.
  \item At a root with $\sigma=0$, by a reduction certificate
    of Section~\ref{sec:perfect}: a chain of reductions that preserve the
    existence of a perfect packing, each confirmed by an exact subset-sum
    test, ending in a residual that is decided by a checked packing, which
    attains the bound $z_0=W/M_1$ of Theorem~\ref{thm:lehmerdual}, or by an
    exhausted search, which proves $z^*=\UB$, the whole chain re-verified by
    independent code.
  \end{thmlist}
  No status of the linear or integer programming solver is a proof.
\end{definition}

The rounded ray $p=\lfloor2^{k_1}y\rfloor\ge0$ of part (iii) is a certificate
when, with the covering rows, the branching rows and the cuts of the node
written as $A'\lambda\ge b'$, the pricing pass gives $p^\top A'_a\le0$ for every
pattern $a$ of the node's space and $p^\top b'>0$: for any feasible $\lambda\ge0$
the chain $0<p^\top b'\le p^\top A'\lambda=\sum_a\lambda_ap^\top A'_a\le0$ would be a
contradiction. Columns bounded by column branching enter the test through
their bounds, every term is an integer, and a ray that fails the test adds
its best pattern as a column, after which column generation resumes. On
the covering rows alone the ray finds the items that the branching rules
and the fixing have left without a covering pattern. A column generation stopped by a
filter, a heuristic or a time limit has not converged, and the rule makes
the correctness of the solver independent of such states.

Column generation at a node must answer one question, prune or branch,
and it can stop as soon as the answer is known. The classical remedy for
tailing off stops once the rounded Lagrangian bound meets the rounded
value of the restricted master~\cite{lubbecke2005}, but at a node whose
relaxation lies just above the target it can lose exactly the unit that
would prune the node. \solver{} uses this rounded stop, with an absolute
tolerance $\tau_r$, only at the root, and the following rule at every
other node.

\begin{definition}[Decision stop]\label{def:decision}
  At a node other than the root, with improving target $T$, column
  generation stops at the first iteration at which one of the following
  holds, in this order.
  \begin{thmlist}
  \item A certified bound $\LB>T+\tau$ is in hand; the node is pruned.
  \item $T-\eta\max(1,T)\le\zrmp\le T$; the node branches on the current
    solution.
  \item Pricing at the fine scale finds no column; the node branches on
    the current solution.
  \item A safeguard fires, a number of iterations in which neither $\zrmp$
    falls nor the certified bound rises, or a time cap; the node branches.
  \end{thmlist}
\end{definition}

The lower end of the window in case (ii) is the strict-below rule. On a
node whose relaxation lies strictly below $T$, the first linear program
over the inherited columns is already below $T$, so case (ii) holds at once
and its branching candidate is an arbitrary vertex of a large face;
branching on it can double the tree at every level. When
\begin{equation}\label{eq:strictbelow}
  \zrmp<T-\eta\max(1,T)
\end{equation}
for the relative tolerance $\eta>0$, column generation therefore
continues to convergence, case (iii), and the node branches on the
converged solution. A node whose column generation stops with an integral
solution has no branching candidate. Its packing, checked bin by bin,
replaces the incumbent if it has fewer bins, after which the node is
closed if the incumbent meets its bound and otherwise returns to the queue
with the lower target; otherwise the node is left unresolved with
its certified bound, and \solver{} reports optimality only when the
incumbent meets the bound of every unresolved node.

\begin{proposition}[Validity of the decision stop]\label{prop:decision}
  At a node with improving target $T$ the decision stop has the following
  properties.
  \begin{thmlist}
  \item Every stop is valid: a prune in case (i) discards no improving
    solution, and the branchings of the other cases partition the integer
    solutions of the node.
  \item Once $\zrmp\le T$, no certified bound of the node over the same
    rows reaches $T+1$, and further iterations cannot prune the node.
  \item If the stop is case (iii), with optimal duals of the restricted
    master and a fine-scale pass that finds no pattern with
    $p^\top a>2^{k_1}(1+\varepsilon)$ for the column acceptance threshold
    $\varepsilon$, and $\zrmp>T$, then
    $|\zlp-T|\le\tau+(T+\tau)(\varepsilon+N_{\max}2^{-k_1})$.
  \end{thmlist}
\end{proposition}

\begin{proof}
  \pitem{i} A prune in case (i) rests on a certified bound $\LB>T+\tau$,
  which is a lower bound on the relaxation of the node by
  Corollary~\ref{cor:bounds}, since $\bar\nu\ge\nu(\pi)$; so the node
  contains no solution of value at most $T$. A branching on $y_{ij}$,
  $z_{ik}$ or a column $\lambda_a$ partitions the integer solutions of the node whatever the
  fractional point at which it is taken, so every branching stop is valid.
  \pitem{ii} As long as the rows of the node do not change, its relaxation
  value $\zlp$ does not depend on the iteration, and once $\zrmp\le T$ it
  satisfies $\zlp\le\zrmp\le T$; every certified bound $\LB$ over these rows
  is at most $\zlp$, so $\lceil\LB-\tau\rceil\le T$.
  \pitem{iii} Optimal duals give $\dv^\top\pi=\zrmp$, and the fine-scale pass
  gives $P^*\le2^{k_1}(1+\varepsilon)$, so that the definition of $\bar\nu$ in
  inequality~\eqref{eq:scaling} at scale $2^{k_1}$ yields
  $[\bar\nu]_+\le\varepsilon+N_{\max}2^{-k_1}\eqdef\varepsilon'$. Since case (i)
  did not occur, part (ii) of Lemma~\ref{lem:resolution} gives
  $\zrmp\le(T+\tau)(1+\varepsilon')$, and part (i) of
  Corollary~\ref{cor:bounds} gives $\zlp\ge\zrmp/(1+\varepsilon')$. Hence
  $\zlp-T\le\zrmp-T\le\tau+(T+\tau)\varepsilon'$, and, as $\zrmp>T$,
  $T-\zlp<T-T/(1+\varepsilon')\le T\varepsilon'$.
\end{proof}

Part (iii) describes the nodes that no bound can settle: when $\zlp=T$,
rounding can hold $\zrmp$ a hair above $T$, and when $\zlp\in(T,T+\tau]$
neither case (i) nor case (ii) can occur. Such nodes are decided by
branching, and the safeguard ends their column generation when it stalls.
Condition~\eqref{eq:precision} says how accurately the pricing value and
the duals must be known to certify the lift $\zrmp-T$ of a node, so the
solver derives the column acceptance threshold and the dual tolerance of
the linear programming solver from that lift, where da Silva and
Schouery~\cite{dasilva2025} use a blanket scaling.

The tree branches on pairs of item types and adds rank-1 cuts, and both
enter the pricing problem. For a solution $\lambda$ and types $i\ne j$ let
$y_{ij}=\sum_a\lambda_a[a_i\ge1][a_j\ge1]$. When $y_{ij}$ is fractional the
children impose $y_{ij}\le\lfloor y_{ij}\rfloor$ and
$y_{ij}\ge\lceil y_{ij}\rceil$, which for unit demands is the branching of
Ryan and Foster~\cite{ryan1981}; for types with $d_i>1$ the tree also
branches on $z_{ik}=\sum_a[a_i\ge k]\lambda_a$, the branching of
Vanderbeck~\cite{vanderbeck2011}. When every $y_{ij}$ and $z_{ik}$ is
integral but $\lambda$ is not, the node branches on a fractional column,
$\lambda_a\le\lfloor\lambda_a\rfloor$ and $\lambda_a\ge\lceil\lambda_a\rceil$, and
pricing in the first child excludes that pattern with a forbidden-aware
bounded knapsack, while the column's upper bound enters the certified
bounds through its dual. Every fractional solution thus has a branching
candidate. In the exact-coverage subtrees of the
$\sigma=0$ classes the pair rows become pricing rules, forbidding a pattern
to contain both items of an ``apart'' pair or exactly one item of a
``together'' pair. A Chv\'atal--Gomory cut
$\sum_a\lfloor u^\top a\rfloor\lambda_a\le\lfloor u^\top\dv\rfloor$ with
multipliers $u\in[0,1)^n$ is valid for every integer solution with exact
coverage. A cut may remove a packing that over-covers an item, but never
all improving ones: deleting surplus copies turns an improving packing
into one with exact coverage and no more bins, the branching
rules partition the integer solutions, so this packing lies in one node of
the tree, and there every cut holds for it. \solver{} separates the
subset-row triples of Jepsen et al.~\cite{jepsen2008} over the fractional
support and rank-1 cuts with four and five types over the types that
co-occur in it. For cuts $G\lambda\le r$ with duals $\gamma\le0$ and
$\rc_{\pi,\gamma}(a)=1-\pi^\top a-\gamma^\top G_a$, the substitution that
proves Theorem~\ref{thm:budget}, with
$\gamma^\top G\lambda=\gamma^\top r-\gamma^\top(r-G\lambda)$, gives
\begin{equation}\label{eq:cutsbudget}
  \one^\top\lambda=\dv^\top\pi+r^\top\gamma+\sum_a\lambda_a\rc_{\pi,\gamma}(a)+\pi^\top e(\lambda)-\gamma^\top(r-G\lambda)
\end{equation}
for every $\lambda\ge0$ with $A\lambda\ge\dv$ and $G\lambda\le r$. Every term
of identity~\eqref{eq:cutsbudget} after the first two is nonnegative when $\rc_{\pi,\gamma}\ge0$, the last
because $\gamma\le0$ and $r-G\lambda\ge0$, so every bound and fixing rule
above holds with $\LB=\dv^\top\pi+r^\top\gamma$, which also takes the place
of $\dv^\top\pi$ in Lemma~\ref{lem:resolution},
Definition~\ref{def:certrule} and Proposition~\ref{prop:decision}. The
classical way to price the rule and cut terms is a dynamic program over
the capacity with one state bit per cut, which pays a factor $2^m$ for $m$
cuts. The following lemma removes that factor.

\begin{lemma}[Type-split pricing]\label{lem:typesplit}
  At a node let the active rules and cuts involve the type set $I$, let
  $\bar I$ be the other types, and let the pattern space of the node be the
  patterns $a\le\bar a$, for the caps $\bar a$ of the node, with loads in
  $[M_1-\sigma,M_1]$. Write a pattern as $a=(a_I,a_{\bar I})$ and its value
  at the node as $v(a)=\pi_I^\top a_I+\varphi(a_I)+\pi_{\bar I}^\top a_{\bar I}$,
  where $\varphi$ collects every rule and cut term and depends on $a_I$
  alone. Then
  \begin{equation}\label{eq:typesplit}
    \max_av(a)=\max_{a_I\le\bar a_I}\Bigl\{\pi_I^\top a_I+\varphi(a_I)+\max_{M_1-\sigma\le\ell\le M_1}\omega\bigl(\ell-\wv_I^\top a_I\bigr)\Bigr\},
  \end{equation}
  where
  $\omega(\ell)=\max\{\pi_{\bar I}^\top a_{\bar I}:\ \wv_{\bar I}^\top a_{\bar I}=\ell,\ a_{\bar I}\le\bar a_{\bar I}\}$,
  with $\omega(\ell)=-\infty$ when no such vector exists, is one bounded
  knapsack over the types of $\bar I$.
\end{lemma}

\begin{proof}
  Every rule and cut involves types of $I$ only, so its contribution to
  the value of $a$ is a function of $a_I$, and the load of $a$ is
  $\wv_I^\top a_I+\wv_{\bar I}^\top a_{\bar I}$. For a fixed $a_I$ the part
  $a_{\bar I}$ ranges over the vectors $a_{\bar I}\le\bar a_{\bar I}$ whose
  load $\ell-\wv_I^\top a_I$ brings the total load $\ell$ into
  $[M_1-\sigma,M_1]$, and the best of them has value
  $\max_\ell\omega(\ell-\wv_I^\top a_I)$; maximizing over $a_I$ gives
  equation~\eqref{eq:typesplit}.
\end{proof}

The outer maximization enumerates the $I$-subpatterns and evaluates
$\varphi$ on each, at a cost linear in the number of rules and cuts, while the inner maximum is read off one bounded knapsack
over $\bar I$, computed once per call. The number of rules and cuts thus
enters the cost linearly instead of as the factor $2^m$, and the
enumeration is governed by the types that the rules and cuts involve,
which are few, since they live where the relaxation is fractional. An
enumeration cut short by a cap returns the larger of the best value found
and the largest upper bound over its unexplored subtrees, so that the value
of the call remains a valid $\bar\nu$ and can also close a node. With the
rule and cut duals scaled and rounded down like the item duals, and their
rounding error added to that of the item profits in
inequality~\eqref{eq:scaling}, the value is exact in integer arithmetic, so one call
at the fine scale is also the certifying pass of
Definition~\ref{def:certrule}, and the $N_{\rm col}$ best $I$-subpatterns
with their completions give $N_{\rm col}$ distinct columns per call.
What cuts can achieve is known exactly.

\begin{proposition}[Cuts at a tight relaxation]\label{prop:nocuts}
  If $z^*=\zlp$, then adding to the relaxation any family of inequalities
  satisfied by one optimal integer solution leaves its value at $\zlp$.
\end{proposition}

\begin{proof}
  The relaxation with the added inequalities is a relaxation of a problem
  whose optimum $z^*$ is attained by a solution satisfying them, so
  $\zlp\le z^{\rm cuts}_{\rm LP}\le z^*=\zlp$.
\end{proof}

On every AI instance $z^*=L_1=\zlp$, so no cut can act there, while on an
ANI instance cuts are the only root device of the relaxation that can
prove the bound; since the status of an instance is unknown, the root cut
loop is budgeted rather than exhaustive.

Algorithm~\ref{alg:lbp} assembles the parts. Its order matters: at a root
with $\sigma=0$ the solver must decide how much time the tools that look
for a perfect packing receive before the column generation that may prove
the bound.

\begin{algorithm}[t]
\caption{Lehmer branch-and-price}\label{alg:lbp}
\small
\begin{algorithmic}[1]
\State compute $L_1$, $L_2$, $M_1$ and the Lehmer bound $z_0=W/M_1$
\State run FFD and BFD; $\UB\gets$ best packing
\If{$\UB=L_1+1$ and $\sigma=0$}
  \State run the tools: reductions, first search slice, backdoor peeling; \Return $\UB$ if decided
\EndIf
\State run the rest of the portfolio until $\UB=\max(L_1,L_2)$; seed the pool with all bins
\If{$\UB=\max(L_1,L_2)$} \Return $\UB$ \EndIf
\If{$\sigma=0$ and the tools have not run}
  \State run them; \Return $\UB$ if decided
\EndIf
\If{the tools ran} replace the instance by the residual of their reductions \EndIf
\State root column generation from the centre $w/M_1$ until $\LB_{\rm node}=\UB$ or a stop rule applies
\State certify the root bound by a fine-scale pricing pass; \Return $\UB$ if $\LB_{\rm node}=\UB$
\State root cut loop; cover, dives and vertex restarts; at $\sigma=0$ the second search slice
\State $\mathcal Q\gets\{\text{root}\}$
\While{$\mathcal Q\ne\emptyset$ and time remains}
  \State take a node of smallest bound from $\mathcal Q$, plunging depth first
  \State column generation with the decision stop, reduced-cost fixing, pricing rules and cuts
  \If{the node is closed by the certification rule} \State \textbf{continue} \EndIf
  \State primal heuristics at the node; update $\UB$ and prune $\mathcal Q$
  \State branch on a fractional pair of types, on $z_{ik}$ or on a column, or leave an integral node unresolved
\EndWhile
\State \Return $\UB$, proved optimal when $\mathcal Q=\emptyset$ and $\UB$ meets the bound of every unresolved node
\end{algorithmic}
\end{algorithm}

First fit decreasing and best fit decreasing run first. When their best
packing has $L_1+1$ bins and leaves $\sigma=0$, as on instances with
$W=L_1c$ and a full bin, the reductions, a first slice of the search and
backdoor peeling run before the rest of the portfolio; this is the path of
every AI and ANI instance. Otherwise the rest of the portfolio runs, and
the tools run after it whenever $\sigma=0$.
A search left undecided by its first slice receives a second one after
the root's column generation, cuts and dives, so that running the tools
first costs little when the relaxation decides the instance. The root
column generation ends when $\LB_{\rm node}$ reaches $\UB$, when the
rounded stop $\lceil\LB_{\rm Lagr}-\tau\rceil\ge\lceil\zrmp-\tau_r\rceil$
applies, or when the relative gap between $\zrmp$ and $\LB_{\rm Lagr}$
falls below a tolerance.

The pricing problem is the bounded knapsack with scaled integer profits,
solved by the COMBO algorithm of Martello, Pisinger and
Toth~\cite{martello1999} on its binary expansion, as in Delorme and
Iori~\cite{delorme2020}; when its optimum has a load below $M_1-\sigma$,
which part (ii) of Theorem~\ref{thm:slack} excludes from every improving
solution, the pricer is the window pricer, an exact-load dynamic program over the
window of loads $[M_1-\sigma,M_1]$ that returns the $N_{\rm col}$ best distinct improving
patterns of a pass. Every improving incumbent of the search is harvested
as a column. The column generation is stabilized by Wentges
smoothing~\cite{wentges1997}, pricing at $\alpha\pi^\circ+(1-\alpha)\pi$
with the centre $\pi^\circ$ the best Lagrangian point so far and $\alpha$
adapted as by Pessoa et al.~\cite{pessoa2018}; the centre starts at the
Lehmer dual $w/M_1$, which is dual feasible, so the Lagrangian bound is
valid from the first iteration. The dual-optimal inequalities of
Val\'erio de Carvalho~\cite{valerio2005} and Ben Amor et
al.~\cite{benamor2006} for pairs and triples of types enter the root
master as exchange columns, fixed to zero at every node with branching
constraints.

The primal heuristics run after the root column generation and at the
nodes. The cover rounds up the linear relaxation of the set multicover
over the pool and removes the over-coverage from the bins of lowest
rigidity, where the rigidity of a pattern $a$ with largest item $j$ is
$\bar\rho_a=\int_1^2\rho_a(s)\dd s$ with
\begin{equation*}
  \rho_a(s)=K_a(s)-\log\bigl(w_j^{s}+(c-w_j)^{s}\bigr)\;\le\;0 .
\end{equation*}
For $s>1$ the inequality is an equality exactly when the other items of
$a$ form one item of weight $c-w_j$, or $a=\{j\}$ with $w_j=c$, so
$\bar\rho_a$ measures how far a bin is from a two-item full bin. A dive
sorts the columns of the cover, or of the support of the relaxation, by
the score $\log\lambda_a+\bar\rho_a$, and three children fix prefixes of
decreasing length; each child re-solves by column generation from the
pool and is pruned as soon as its fixed columns exceed the reduced-cost
budget of part (iii) of Corollary~\ref{cor:bounds} or their total slack
exceeds $\sigma$. At $\sigma=0$ the relaxation lies on a massively degenerate face,
and whether a dive succeeds depends on the vertex that the simplex method
returned; re-solving the restricted master after a small random perturbation of
its costs typically returns another vertex of the face, and a few dives
from such vertices follow a failed root dive.

The tree is best-bound with depth-first plunging. A node is pruned by its
bound~\eqref{eq:lbnode} before any column generation, its column
generation runs with the decision stop and reduced-cost fixing, and in the
exact-coverage subtrees the pair rows become the pricing rules of
Lemma~\ref{lem:typesplit}. On instances with the integer round-up
property the root bound almost always equals the optimum, and a tree is a
search for the packing that the primal heuristics conduct; on instances
without it every node must be pruned, which is where the decision stop,
the cuts and the strict-below rule act.
The parts assembled above make the optimality reported by \solver{} a
theorem.

\begin{theorem}[Soundness]\label{thm:soundness}
  If \solver{} reports $\UB$ as optimal, then $\UB=z^*$.
\end{theorem}

\begin{proof}
  $\UB$ is the number of bins of a packing checked bin by bin, so
  $z^*\le\UB$. Suppose that a packing with fewer than $\UB$ bins exists, and
  let $\lambda^\circ$ be the integer solution with exact coverage and no more
  bins obtained from it by deleting surplus copies. If the tools replaced
  the instance by a residual, they did so at a root with $\sigma=0$, where
  every improving solution is a perfect packing for the capacity $M_1$ by
  part (iii) of Theorem~\ref{thm:slack}, through the reductions of
  Lemmas~\ref{lem:exchange} and~\ref{lem:forced}, which preserve the
  existence of such a packing; so we may take $\lambda^\circ$ in the
  residual. The branching rules partition the integer solutions, so
  $\lambda^\circ$ lies in a node $\nu$ that was closed by
  Definition~\ref{def:certrule} or left unresolved: the root if \solver{}
  returned before the tree, and otherwise a leaf, since the run ended with
  an empty queue. Let $f$ be the number of fixed bins of $\nu$ and $t$ the
  number of the other bins of $\lambda^\circ$, so that $f+t<\UB$ for the
  reported $\UB$, and $t\le T$ for the improving target $T$ of $\nu$ at any
  time, since the incumbent never increases. At $\nu$, $\lambda^\circ$
  satisfies the branching rows, the cuts, which hold for solutions with
  exact coverage, the reduced-cost fixing, by part (iii) of
  Corollary~\ref{cor:bounds}, since $\lambda^\circ$ has no over-coverage, and
  the window of loads, by part (ii) of Theorem~\ref{thm:slack}. So it is
  feasible for the relaxation of $\nu$, and $\zlp\le t$. Every bound of
  Corollary~\ref{cor:bounds} at $\nu$ is therefore at most $t$, the
  Lagrangian bound because $U\ge T+1>\zlp$, and since the bound on its
  rounding error is subtracted, its computed value $\LB$ has
  $\lceil\LB-\tau\rceil\le t$. Moreover $L_1$ and $L_2$ are at most $f+t$, and
  so is every bound inherited from an ancestor, whose solutions include
  $\lambda^\circ$. Hence $\LB_{\rm node}\le f+t<\UB$. We check each outcome.
  \pitem{i} A packing of $\nu$ attaining $\LB_{\rm node}$ has at most $f+t$
  bins, so the reported $\UB$ would not exceed $f+t$.
  \pitem{ii} $\LB_{\rm node}<\UB$, and $\lceil\LB-\tau\rceil\le t\le T$.
  \pitem{iii} $\sigma\ge0$ by part (i) of Theorem~\ref{thm:slack}, every item
  is covered by a pattern of $\lambda^\circ$, and $\lambda^\circ$, feasible for
  the master of $\nu$, rules out a certified Farkas ray by the chain that
  certifies it.
  \pitem{iv} This closure happens only at the root, where every improving
  solution is a perfect packing for the capacity $M_1$. A checked packing
  found there has $z_0$ bins, and $z_0\le z^*$ by Theorem~\ref{thm:lehmerdual};
  the search branches on a blocker set by Lemma~\ref{lem:blocker} or on all
  completions of an item, so an exhausted search contradicts the existence
  of $\lambda^\circ$.

  \noindent If $\nu$ was left unresolved, its bound is at most $f+t<\UB$, and
  optimality is not reported. Hence $\nu$ could neither be closed nor left
  unresolved, a contradiction.
\end{proof}

\section{Computational Results}\label{sec:results}

This section compares \solver{} with the general-purpose methods that
lead on the standard benchmark, and then separates the contributions of
the perfect-packing tools, the entropic tools and the components of the
engine.

The benchmark is the set of $6\,195$ instances of BPPLIB~\cite{bpplib}:
the $5\,955$ instances of the review of Delorme, Iori and
Martello~\cite{delorme2016} and the $240$ GI instances of Gschwind and
Irnich~\cite{gschwind2016}, with capacities up to $1.5\cdot10^6$. The
classical classes are Falkenauer~U and T~\cite{falkenauer1996}, Scholl~1,
2 and 3~\cite{scholl1997}, Schwerin~1 and 2~\cite{schwerin1997},
W\"ascher~\cite{wascher1996}, Hard28~\cite{schoenfield2002}, Random and
GI. The AI and ANI classes are ten groups of fifty instances with $201$ to
$1\,003$ items and capacities up to $78\,080$, whose optima are known by
construction, $L_1$ on AI and $L_1+1$ on ANI.

We compare with published results of NF-F, the reduced network flow method
of de Lima, Iori and Miyazawa~\cite{delima2023} (Xeon E3-1245~v5,
$3.5$~GHz), and dSS, the branch-cut-and-price of da Silva and
Schouery~\cite{dasilva2025} (Xeon E5-2630~v4, $2.2$~GHz), both with
Gurobi and a one-hour limit, and of BELOV, the branch-and-cut-and-price of
Belov and Scheithauer~\cite{belov2006}, the strongest code of the
review~\cite{delorme2016,delorme2020}. No other code in the literature
solves more instances on these
classes~\cite{wei2020,delorme2020,delima2023,delorme2026}. \solver{} is
written in C++. Its pricing uses a port of the COMBO code of
Martello, Pisinger and Toth~\cite{martello1999}, the Lehmer duals and
heuristics a port of the bounded knapsack code of
Pisinger~\cite{pisinger2000}, its subset-sum and exact-cover kernels are
its own, and it uses the open-source solver
HiGHS~1.11~\cite{huangfu2018} for every linear and integer program. It
runs on one thread of an Intel Xeon w5-2465X with a limit of $600$
seconds per instance; times are reported as measured, each method on its
own machine. A run counts as solved only when its closure satisfies
Definition~\ref{def:certrule}. Every optimum agrees with the known one,
every reduction certificate passes the independent checker, and about
$1.3\cdot10^5$ numerical checks of the identities, bounds and stopping rules
of Sections~\ref{sec:lehmer} to~\ref{sec:bp} found no violation beyond
floating-point rounding. The review ran BELOV only on the instances whose
$L_2$ bound and best fit decreasing value differ; we count the others as
solved in zero time and weight overall times by class size. Bold marks
the best general-purpose result, ours when tied. The parameters of
\solver, the same on every instance, are listed in Table~\ref{tab:params}.
The perfect-packing tools were developed on the AI and ANI instances, so
on these classes the like-for-like comparison with the general-purpose
methods is the engine without them.

\begin{table}[t]
\caption{Parameters of \solver.}\label{tab:params}
\centering\tablesmall
\begin{tabular}{@{}l>{\raggedright\arraybackslash}p{0.62\textwidth}@{}}
\toprule
value & role \\
\midrule
$\tau=10^{-9}$, $\eta=10^{-6}$ & certification tolerance; relative tolerance of the strict-below rule~\eqref{eq:strictbelow} \\
$k_0=30$, $k_1=52$ & coarse pricing scale; fine scale of the certifying pass, with $128$-bit products \\
$\min(10^{-9},0.1\,\delta/T)$ & tightened tolerances: column acceptance at nodes with lift $\delta=\zrmp-T<10^{-5}$, with linear programming tolerances $10^{-10}$ and the objective scaled by $2^8$ \\
$\alpha\in[0,0.9]$ & range of the Wentges weight \\
$\tau_r=10^{-3}$, $10^{-6}$ & absolute tolerance of the rounded stop and relative gap $(\zrmp-\LB_{\rm Lagr})/\max(1,\zrmp)$ that end the root column generation \\
$N_{\rm col}=30$ & distinct improving patterns returned per pass of the window and type-split pricers \\
$c>10^5$ & capacities at which the subset-sum heuristics are skipped \\
$15$, $30\%$, $4$ & root cut loop: rounds, share of the remaining time, rounds without lift \\
$1$, $0.75$, $0.5$; $35\%$ & fractions of the sorted columns fixed by the children of a dive; work budget of the dives at $\sigma=0$ nodes \\
$2$, $10^{-6}$ & vertex restarts after a failed root dive, and their relative cost perturbation \\
\midrule
$0.5$~s, $16.9$~s, $10^6$ & search of Section~\ref{sec:perfect}: first slice, total time and node limit per root \\
$\min(30\text{ s},5\%)$ & budget of backdoor peeling, in seconds or as a share of the time limit \\
$\kappa=60$; $\mu=12$, $4$, $150$ & residual size at which the reserved cascade stops; conflicting completions per triple, items per candidate reserve, candidates per stall \\
$b=5$, $n_V=8$; $5+5$ & beam width and largest reserve; candidates rolled out per stall, the best five and five further single items \\
$4$; $5$~s, $10^6$ & exempt large weights after which a reserve is rejected; search on each residual of backdoor peeling \\
\bottomrule
\end{tabular}
\end{table}

On the classical classes other than GI (Table~\ref{tab:solved}), NF-F,
dSS and \solver{} solve all $5\,455$ instances, and BELOV all but $23$
Falkenauer~T instances within the one-minute limit of the review.
\solver{} has the lowest mean time on every class except Falkenauer~T,
where dSS is faster, and Falkenauer~U, where BELOV is, and over these
classes as a whole, $0.06$ seconds against $0.08$ for dSS. On GI
(Table~\ref{tab:gi}, with the results of dSS and of BELOV as run by
Delorme and Iori~\cite{delorme2020}, both reported by da Silva and
Schouery~\cite{dasilva2025}), \solver{} and BELOV solve all $240$
instances, and \solver{} has the lowest mean time over the class.

\begin{table}[t]
\caption{Classical classes other than GI: instances solved and mean time in seconds.}\label{tab:solved}
\centering\tablesmall
\setlength{\tabcolsep}{3.6pt}
\begin{tabular}{@{}lrrrrrrrrr@{}}
\toprule
 & & \multicolumn{2}{c}{NF-F} & \multicolumn{2}{c}{dSS} & \multicolumn{2}{c}{BELOV} & \multicolumn{2}{c}{\solver{}} \\
\cmidrule(lr){3-4}\cmidrule(lr){5-6}\cmidrule(lr){7-8}\cmidrule(lr){9-10}
class & $n$ & solved & time & solved & time & solved & time & solved & time \\
\midrule
Falkenauer T & 80 & 80 & 0.30 & 80 & \textbf{0.10} & 57 & 24.70 & \textbf{80} & 0.44 \\
Falkenauer U & 80 & 80 & 0.10 & 80 & 0.01 & 80 & \textbf{0.00} & \textbf{80} & 0.01 \\
Hard28 & 28 & 28 & 23.60 & 28 & 6.40 & 28 & 7.30 & \textbf{28} & \textbf{3.85} \\
Random & 3\,840 & 3\,840 & 0.90 & 3\,840 & 0.05 & 3\,840 & 0.15 & \textbf{3\,840} & \textbf{0.05} \\
Scholl & 1\,210 & 1\,210 & 1.40 & 1\,210 & 0.05 & 1\,210 & 0.18 & \textbf{1\,210} & \textbf{0.01} \\
Schwerin & 200 & 200 & 0.20 & 200 & 0.02 & 200 & 1.20 & \textbf{200} & \textbf{0.00} \\
W\"ascher & 17 & 17 & 161.20 & 17 & 0.30 & 17 & 0.10 & \textbf{17} & \textbf{0.03} \\
\midrule
all & 5\,455 & 5\,455 & 1.58 & 5\,455 & 0.08 & 5\,432 & 0.59 & \textbf{5\,455} & \textbf{0.06} \\
\bottomrule
\end{tabular}
\end{table}

\begin{table}[t]
\caption{GI instances: instances solved and mean time in seconds.}\label{tab:gi}
\centering\tablesmall
\setlength{\tabcolsep}{4pt}
\begin{tabular}{@{}lrrrrrrr@{}}
\toprule
 & & \multicolumn{2}{c}{dSS} & \multicolumn{2}{c}{BELOV} & \multicolumn{2}{c}{\solver{}} \\
\cmidrule(lr){3-4}\cmidrule(lr){5-6}\cmidrule(lr){7-8}
subclass & $n$ & solved & time & solved & time & solved & time \\
\midrule
AA & 60 & 60 & 14.50 & 60 & \textbf{2.80} & \textbf{60} & 3.69 \\
AB & 60 & 60 & 7.90 & 60 & 10.90 & \textbf{60} & \textbf{6.51} \\
BA & 60 & 59 & 102.00 & 60 & \textbf{2.80} & \textbf{60} & 6.15 \\
BB & 60 & 60 & 25.70 & 60 & 10.50 & \textbf{60} & \textbf{10.14} \\
\midrule
all & 240 & 239 & 37.50 & 240 & 6.75 & \textbf{240} & \textbf{6.62} \\
\bottomrule
\end{tabular}
\end{table}

Table~\ref{tab:classes} gives the median and maximum time in seconds, the
fraction of instances closed at the root, the mean number of tree nodes
and the peak memory in megabytes. The root closes $98.9\%$ of the classical
instances: the portfolio meets $\max(L_1,L_2)$ on $2\,879$, and of the other
roots $1\,277$ close by the certified bound, $1\,407$ by a packing of the
primal heuristics and $68$ by the perfect-packing tools. The root bound
equals the optimum on all but nine instances, five of Hard28, two of
W\"ascher and two of Random, whose trees of $3$ to $1\,127$ nodes close their
leaves by certified bounds; every other tree is a search for a packing,
which is why trees are large only on Hard28. The portfolio reaches the
optimum on $4\,169$ of the $5\,455$ instances outside GI, and no member
dominates: on the $2\,044$ Random instances on which every member runs,
SS-1.5 is optimal on $965$, first fit decreasing on $639$ and the portfolio
on $1\,041$. The Lehmer dual is weak as a bound, $\lceil z_0\rceil$ never
exceeding $L_2$, but close to optimal as a dual solution, with $z_0/\zlp$
averaging $0.997$ over the $633$ roots at which column generation
converges.

\begin{table}[t]
\caption{Profile of \solver{} on the classical classes.}\label{tab:classes}
\centering\tablesmall
\setlength{\tabcolsep}{4pt}
\begin{tabular}{@{}lrrrrrr@{}}
\toprule
class & $n$ & median & max & root & nodes & MB \\
\midrule
Falkenauer T & 80 & 0.26 & 3.43 & 0.99 & 1.61 & 74 \\
Falkenauer U & 80 & 0.00 & 0.04 & 1.00 & 0.00 & 15 \\
Hard28 & 28 & 0.99 & 38.00 & 0.21 & 110.71 & 22 \\
Random & 3\,840 & 0.01 & 37.29 & 0.99 & 0.16 & 35 \\
Scholl 1 & 720 & 0.00 & 0.03 & 1.00 & 0.00 & 15 \\
Scholl 2 & 480 & 0.00 & 0.28 & 1.00 & 0.01 & 19 \\
Scholl 3 & 10 & 0.80 & 1.32 & 1.00 & 0.00 & 19 \\
Schwerin 1 & 100 & 0.00 & 0.02 & 1.00 & 0.00 & 13 \\
Schwerin 2 & 100 & 0.00 & 0.02 & 1.00 & 0.00 & 13 \\
W\"ascher & 17 & 0.00 & 0.12 & 0.82 & 1.35 & 13 \\
GI & 240 & 2.75 & 315.91 & 0.95 & 2.04 & 210 \\
\midrule
all & 5\,695 & 0.00 & 315.91 & 0.99 & 0.77 & 210 \\
\bottomrule
\end{tabular}
\end{table}

The augmented classes are where the methods differ
(Table~\ref{tab:aiani}, whose groups are named by the number of items of
their ANI instances, the AI instances having one item more; mean times count unsolved runs at the limit, and the
column of the engine, \solver{} without the perfect-packing tools, is not
set in bold). \solver{} solves all $500$ instances, in $0.17$ seconds on
average and $5.6$ at most, while NF-F and dSS solve $473$ and $487$ within
one hour. The engine alone solves $430$ within ten minutes, the figure to
set against the general-purpose codes. The class-specific algorithms of
da Silva et al.~\cite{dasilva2026}, which exploit the construction of the
two classes, also solve all $500$ within one hour, with the same mean time
as \solver{} over AI, $0.27$ seconds, and $2.69$ seconds against $0.06$
over ANI.

\begin{table}[t]
\caption{AI and ANI classes: instances solved out of 50 and mean time in seconds.}\label{tab:aiani}
\centering\tablesmall
\setlength{\tabcolsep}{4pt}
\begin{tabular}{@{}llrrrrrrrr@{}}
\toprule
 & & \multicolumn{2}{c}{NF-F} & \multicolumn{2}{c}{dSS} & \multicolumn{2}{c}{\solver{} engine} & \multicolumn{2}{c}{\solver{}} \\
\cmidrule(lr){3-4}\cmidrule(lr){5-6}\cmidrule(lr){7-8}\cmidrule(lr){9-10}
class & group & solved & time & solved & time & solved & time & solved & time \\
\midrule
ANI & 201 & 50 & 3.00 & 50 & 0.40 & 50 & 1.58 & \textbf{50} & \textbf{0.00} \\
 & 402 & 50 & 24.90 & 50 & 2.20 & 50 & 3.30 & \textbf{50} & \textbf{0.01} \\
 & 600 & 50 & 140.70 & 50 & 11.80 & 46 & 63.49 & \textbf{50} & \textbf{0.04} \\
 & 801 & 49 & 393.20 & 50 & 57.00 & 48 & 55.57 & \textbf{50} & \textbf{0.07} \\
 & 1002 & 43 & 1302.50 & 47 & 374.00 & 37 & 249.62 & \textbf{50} & \textbf{0.16} \\
 & all & 242 & 372.90 & 247 & 89.10 & 231 & 74.71 & \textbf{250} & \textbf{0.06} \\
\midrule
AI & 201 & 50 & 2.00 & 50 & 0.40 & 50 & 0.68 & \textbf{50} & \textbf{0.17} \\
 & 402 & 50 & 25.20 & 50 & 5.00 & 44 & 76.47 & \textbf{50} & \textbf{0.15} \\
 & 600 & 49 & 192.40 & 50 & 57.10 & 35 & 192.43 & \textbf{50} & \textbf{0.31} \\
 & 801 & 46 & 566.50 & 48 & 223.70 & 38 & 162.33 & \textbf{50} & \textbf{0.30} \\
 & 1002 & 36 & 1577.10 & 42 & 794.90 & 32 & 258.43 & \textbf{50} & \textbf{0.44} \\
 & all & 231 & 472.60 & 240 & 216.20 & 199 & 138.07 & \textbf{250} & \textbf{0.27} \\
\bottomrule
\end{tabular}
\end{table}

Every AI and ANI instance is closed at the root, and
Table~\ref{tab:aidetail} gives the closures, by an item without a
completion (dead), by the reductions and the search (search), by backdoor
peeling (backdoor) and by the reductions alone (empty), with the seconds
spent in each tool. Each ANI instance is refuted, $45$ by a dead item and
$205$ by a search of at most eight nodes on a residual kernel. On AI, $179$
instances are packed by the reductions and the search, $69$ by backdoor
peeling with reserves of one to four items and two by the reductions
alone, and every packing is checked bin by bin on the original items.

\begin{table}[t]
\caption{Profile of \solver{} on the AI and ANI groups.}\label{tab:aidetail}
\centering\tablesmall
\setlength{\tabcolsep}{2.6pt}
\begin{tabular}{@{}lrrrrrrrrrr@{}}
\toprule
 & \multicolumn{3}{c}{time} & \multicolumn{4}{c}{closures} & \multicolumn{3}{c}{seconds in} \\
\cmidrule(lr){2-4}\cmidrule(lr){5-8}\cmidrule(lr){9-11}
group & mean & median & max & dead & search & backdoor & empty & reductions & search & backdoor \\
\midrule
AI-201 & 0.17 & 0.03 & 0.53 & 0 & 35 & 14 & 1 & 0.00 & 0.16 & 0.00 \\
AI-402 & 0.15 & 0.02 & 0.56 & 0 & 39 & 11 & 0 & 0.01 & 0.13 & 0.01 \\
AI-600 & 0.31 & 0.04 & 5.60 & 0 & 36 & 14 & 0 & 0.02 & 0.16 & 0.12 \\
AI-801 & 0.30 & 0.08 & 0.96 & 0 & 32 & 17 & 1 & 0.05 & 0.18 & 0.05 \\
AI-1002 & 0.44 & 0.16 & 2.49 & 0 & 37 & 13 & 0 & 0.11 & 0.16 & 0.14 \\
ANI-201 & 0.00 & 0.00 & 0.01 & 19 & 31 & 0 & 0 & 0.00 & 0.00 & 0.00 \\
ANI-402 & 0.01 & 0.01 & 0.02 & 6 & 44 & 0 & 0 & 0.01 & 0.00 & 0.00 \\
ANI-600 & 0.04 & 0.03 & 0.04 & 6 & 44 & 0 & 0 & 0.02 & 0.00 & 0.00 \\
ANI-801 & 0.07 & 0.07 & 0.09 & 9 & 41 & 0 & 0 & 0.05 & 0.01 & 0.00 \\
ANI-1002 & 0.16 & 0.15 & 0.28 & 5 & 45 & 0 & 0 & 0.12 & 0.02 & 0.00 \\
\midrule
all & 0.17 & 0.07 & 5.60 & 45 & 384 & 69 & 2 & 0.04 & 0.08 & 0.03 \\
\bottomrule
\end{tabular}
\end{table}

\begin{figure}[t]
\centering
\includegraphics[width=\textwidth]{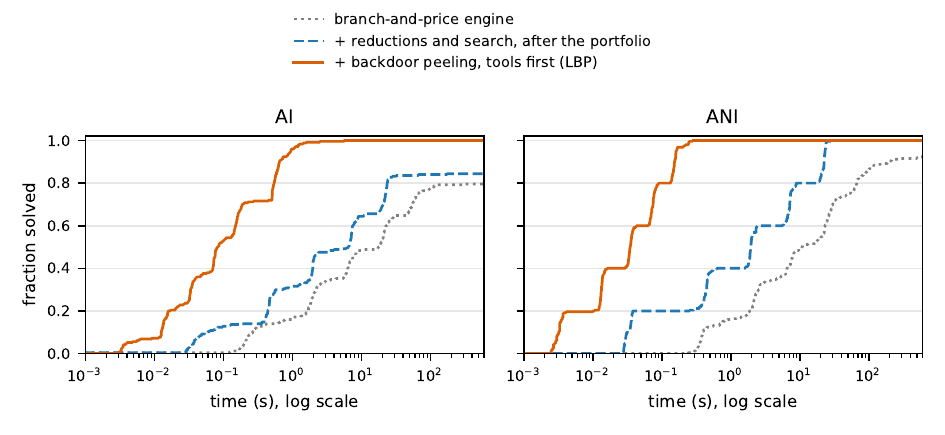}
\caption{Fraction of the AI and ANI instances solved against time by the configurations in the upper part of Table~\ref{tab:machinery}.}\label{fig:aiani}
\end{figure}

Table~\ref{tab:machinery} and Figure~\ref{fig:aiani} separate the
contributions of the tools: the upper part of the table adds them in
turn, and the lower part counts the roots that they decide alone, without
column generation. Run after the whole portfolio, the reductions and the
search close every ANI instance but only $12$ AI instances more than the
engine. Backdoor peeling closes the remaining $39$ AI instances, $37$ of
which that configuration leaves open even in one hour, and running the
tools right after first fit decreasing and best fit decreasing saves most
of the time of that configuration, which spends $99\%$ of its time on ANI
in the portfolio before the tools run. Alone, the reductions with an exact cover of
the residual decide every ANI root and half of the AI roots;
fewest-completion branching decides $41$ more, blocker branching $19$ and
backdoor peeling the remaining $64$.

\begin{table}[t]
\caption{Contribution of the perfect-packing tools on the AI and ANI classes.}\label{tab:machinery}
\centering\tablesmall
\setlength{\tabcolsep}{4pt}
\begin{tabular}{@{}lrrrr@{}}
\toprule
configuration & AI & ANI & total & mean s \\
\midrule
branch-and-price engine & 199 & 231 & 430 & 106.39 \\
+ reductions and search & 211 & \textbf{250} & 461 & 53.32 \\
+ backdoor peeling, tools first (\solver{}) & \textbf{250} & \textbf{250} & \textbf{500} & \textbf{0.17} \\
\midrule
reductions, exact cover of the residual & 126 & \textbf{250} & 376 & \\
+ search, fewest-completion branching & 167 & \textbf{250} & 417 & \\
+ blocker branching & 186 & \textbf{250} & 436 & \\
+ backdoor peeling & \textbf{250} & \textbf{250} & \textbf{500} & \\
\bottomrule
\end{tabular}
\end{table}

The tools succeed because of how the augmented classes are built. Each
ANI instance is a non-IRUP kernel of fifteen items to which triplets that
fill a bin are added one at a time, each large item chosen so that no
subset of the items present before it completes
it~\cite{delorme2016,dasilva2026}. The last triplet therefore satisfies
Lemma~\ref{lem:exchange}, and so does each earlier one once the later ones
are fixed, so the reductions peel the instance down to its kernel unless
an item without a completion closes it first. Each
AI instance is derived from an ANI instance by splitting one item into
two pieces~\cite{dasilva2026}, which give the large items completions that the construction
excluded. On the AI instances that the reductions and the search leave
undecided, the cascade therefore stalls at its first step, fixing a median
of $0.4\%$ of the large items; with the two pieces held back from the
bins of the large items, the reductions peel every AI instance to a
residual of $16$ to $22$ items that packs. The solver does not know the
pieces: on the $248$ AI roots that backdoor peeling packs when run alone,
the reserve it learns is the pair on $57$, one piece on $42$, empty on $131$ and another set on the other
$18$, since the construction admits several reserves that unlock a
packing. On the classical classes the tools act only at the $75$ roots
with $\sigma=0$, mostly of Falkenauer~T, and change the total time by
at most half a percent.

The entropic tools of Section~\ref{sec:entropic} are not part of the
default configuration. A pass of the sum--product knapsack costs a median
of $0.34$ exact pricing calls on the classical roots and $0.77$ on the AI
and ANI roots, and is expensive only at the capacities of GI. Over
$54\,118$ evaluations at root duals, the Gibbs identity and the entropic
bounds held up to one rounding error of $3\cdot10^{-8}$. At the duals at
which the root column generation ends, the effective number
$e^{H(P_{1000})}$ of nearly optimal patterns has median $1\,010$ on the
classical classes against $1.1\cdot10^{15}$ on AI and ANI; since this count
also grows with the capacity, it indicates rather than measures the
degeneracy that makes the augmented classes hard. As root arms on the AI
and ANI instances with the tree disabled (Table~\ref{tab:ablation}, which
reports the roots closed and the mean root iterations and seconds), the
entropic seed, fifty steps of the ascent of Theorem~\ref{thm:smoothed}
followed by $2\,000$ samples of the ensemble, closes the most roots, $339$,
and halves the iterations at a cost of $27$ seconds; the coarse bound, part
(iv) of Theorem~\ref{thm:entbounds} on the coarsened instance, raises the
lower bound earlier, so that the root column generation stops after $6\%$
fewer iterations, and the dives and the cut loop that follow close $8$
more roots, the bound itself closing none, since no lower bound can exceed
$\zlp=T$ at these roots; and smoothing at every iteration fails, reaching
the time limit on most roots and closing only $43$. On a stratified sample of $300$ classical
instances at $60$ seconds, the seed loses $29$ instances, $28$ of them on GI.
The ensemble is therefore a source of bounds, diagnostics and a seed, not
a pricing engine.

\begin{table}[t]
\caption{Ablation of the engine on the AI and ANI classes: root arms (top) and tree arms (bottom).}\label{tab:ablation}
\centering\tablesmall
\setlength{\tabcolsep}{4pt}
\begin{tabular}{@{}lrrr@{}}
\toprule
root arm & roots & iterations & seconds \\
\midrule
engine default & 325 & 1\,671.39 & 50.28 \\
no stabilization & 209 & 14\,814.46 & 290.52 \\
zero initial centre & 328 & 1\,671.86 & 50.44 \\
no exchange columns & 327 & 1\,671.39 & 50.33 \\
no harvesting of incumbents & 325 & 1\,671.39 & 50.35 \\
diverse graph pricer & 337 & 236.82 & 54.42 \\
entropic seed & 339 & 883.44 & 77.25 \\
entropic coarse bound & 333 & 1\,564.21 & 56.34 \\
smoothing at every iteration & 43 & 626.93 & 517.97 \\
portfolio of FFD and BFD & 315 & 2\,203.44 & 72.55 \\
portfolio with SS-2 & 327 & 1\,669.66 & 55.41 \\
\bottomrule
\end{tabular}

\medskip
\begin{tabular}{@{}lrrrr@{}}
\toprule
tree arm & AI & ANI & all & time ratio \\
\midrule
engine default & 27 & 29 & 56 & 1.00 \\
no cuts & 26 & 23 & 49 & 0.96 \\
cuts at the root only & 25 & 27 & 52 & 1.06 \\
subset-row triples only & 26 & 27 & 53 & 1.09 \\
convergence at every node & 26 & 22 & 48 & 1.04 \\
target stop & 27 & 22 & 49 & 1.01 \\
no tightened tolerances & 26 & 25 & 51 & 1.02 \\
no dives & 24 & 28 & 52 & 1.05 \\
column branching first & 26 & 17 & 43 & 1.05 \\
window pricer & 25 & 19 & 44 & 1.25 \\
graph pricer & 25 & 21 & 46 & 1.29 \\
two threads & 27 & 28 & 55 & 0.88 \\
\bottomrule
\end{tabular}
\end{table}

The ablation changes one component of the engine at a time against a
default. Its default
differs from the engine of Tables~\ref{tab:solved} to~\ref{tab:machinery}
in three settings, one column per call of the window pricer, cut rounds at
every node and a time cap on the column generation of a node, and on the
classical sample the arms other than the dives and the portfolio run
without the strict-below rule. On that sample, where the default solves
$299$ instances, the dives are decisive: without them $48$ instances are
lost, $27$ on GI and $12$ on Hard28, where the root bound is tight and the
tree is a search for the packing. Stabilization is worth five instances
and pair branching four, and no other component outside the entropic
tools more than three.

On AI and ANI the root arms run on all $500$ instances and the tree arms
on a stratified sample of $60$, six per group, with the instances solved
per class and the ratio of the geometric means of the times plus one
second over the instances solved by both the arm and the default.
Stabilization dominates the root: without it the root takes nine times
the iterations and closes $116$ fewer roots. The devices that act on the vertex
from which the first dive starts, the entropic seed, a zero initial centre
and a diverse pricer that returns up to thirty columns far apart in the
pattern space from one longest-path pass over the graph whose paths are
the full patterns, close three to fourteen more roots. Cuts are worth six ANI
instances; on AI, where by Proposition~\ref{prop:nocuts} they cannot raise
the bound, removing them loses one instance. The decision stop is worth
seven ANI instances: the two arms that replace it, which end the column
generation of a node at the relative gap of the root, the second also as
soon as $\zrmp\le T$, lose the same seven, and the first also one AI
instance. The relaxation of a prunable node lies of the order
of $10^{-5}$ above the target, and only a stop that continues until the
prune is certified at the fine scale finds it; the tightened tolerances
that let condition~\eqref{eq:precision} hold at such nodes are worth four
ANI instances and one AI instance. Pair branching, whose rules type-split
pricing exploits, is worth twelve ANI instances and one AI instance
against column branching first. Type-split pricing is in turn worth ten
and twelve instances
against the graph pricer, which carries the rules and cuts as states along
the graph of full patterns, and the window pricer, which handles the rules
by a branch-and-bound over window dynamic programs; their calls cost up to
a hundred times and four to ten times as much. Over the $341$ runs that
price with rules or cuts, a call involves $3.8$ types in the median of the
per-run means. The dives are worth three AI instances and one ANI instance,
and a second thread for the primal side shortens the times by $12\%$ but
solves one instance fewer.

\section{Conclusion}\label{sec:conclusion}

Lehmer branch-and-price rests on one observation: read through the Lehmer
transform, a bin is a smooth increasing curve rather than a finite
multiset, and the discrete problem acquires a structure on which the tools
of analysis can act. The curve yields an explicit dual solution, a
monotone representation of the dominance order, a family of heuristics
with first fit decreasing as its limit, and an integer slack identity that
tells when only a perfect packing can improve the best solution known.
Carried to patterns, it gives entropic column generation, with two-sided
bounds at every iteration without exact pricing, a smoothed dual maximized at an explicit rate
and a duality gap that is an entropy. The branch-and-price built on these
foundations closes every node by a certificate and decides the nodes on
which the relaxation is silent by exchange reductions, blocker branching
and backdoor peeling. It solves every instance of the standard
benchmark within ten minutes, the augmented ones in a fraction of a second
on average, and its engine alone solves every classical instance with the
lowest mean time among the methods compared.

A first direction for future work tests the reach of the perfect-packing
tools. The augmented classes are built from chains of exchanges, which is why
they yield to exchange reasoning once a handful of items are held back.
Perfect-packing instances built differently, and fresh instances of the
same construction at larger sizes, are the natural test of their
generality, and the benchmark cannot provide it.

A second direction concerns the transform itself. The cone of
nondecreasing superadditive functions characterizes the fractional
relaxation of the dominance order, and a characterization of dominance
itself remains open; it would tell when one pattern can replace another
in every solution, and so restrict the pattern space before any pricing.
The ensemble, in turn, prices every pattern at once, and pricing rules
that use it to choose diverse columns, or to steer a restricted master
away from its degenerate vertices, could turn the degeneracy that makes
the hardest instances hard from an obstacle into a guide.

\section*{Data Availability}

The source code of Lehmer branch-and-price (version 1.0.1) and of its
certificate checker, the per-instance results and the scripts that
reproduce the tables and the figure are at
\url{https://github.com/chimera1001/lehmer-bp}; the instances are in
BPPLIB~\cite{bpplib}.

\clearpage
\begin{appendices}
\section{Proofs}\label{app:proofs}

This appendix proves the results of the body whose proofs are long or
computational, in the order in which they appear.

\begin{proof}[Proof of Theorem~\ref{thm:ssffd}]
  Distinct types have distinct weights, so a pattern is determined by the
  sequence of its weights in nonincreasing order,
  \begin{equation*}
    a_{(1)}\;\ge\;a_{(2)}\;\ge\;\dots\;\ge\;a_{(D)}\;\ge\;0,
  \end{equation*}
  padded with zeros to length $D$; every pattern of the residual has at
  most $D$ items. We show in turn that the maximizer is the
  lexicographically largest pattern, that this pattern is built greedily,
  and that the greedy bins are those of first fit decreasing.
  \pitem{a} Let $a\ne b$ be patterns whose sequences first differ at
  position $t$, with
  \begin{equation*}
    u=a_{(t)}\;>\;v=b_{(t)}\;\ge\;0,
  \end{equation*}
  so that $a$ is lexicographically larger. All positions before $t$ agree
  and $b_{(r)}\le v$ for $r\ge t$, hence
  \begin{equation*}
    Z_a(s)-Z_b(s)\;\ge\;u^s-\sum_{r\ge t}b_{(r)}^s\;\ge\;u^s-Dv^s ,
  \end{equation*}
  which is positive when $v=0$. For $v>0$ the weights are integers with
  $v\le c$, so $u/v\ge1+1/c$, and since $\ln(1+1/c)\ge1/(c+1)$,
  \begin{equation*}
    \Bigl(\frac uv\Bigr)^s\;\ge\;\Bigl(1+\frac1c\Bigr)^s\;\ge\;e^{s/(c+1)}\;>\;D
    \qquad\text{for }s>(c+1)\ln D .
  \end{equation*}
  Then $u^s>Dv^s$ and $Z_a(s)>Z_b(s)$. For such $s$ the order of the
  values $Z_a(s)$ therefore refines the lexicographic order, and the
  maximizer is the lexicographically maximal pattern.
  \pitem{b} The greedy bin $g$ receives its items in nonincreasing order,
  since the remaining capacity only decreases. Let $a$ be the
  lexicographically maximal pattern, and suppose that $a\ne g$, the
  sequences first differing at position $t$. The item $g_{(t)}$ is the
  largest residual item, apart from $g_{(1)},\dots,g_{(t-1)}$, that fits in
  the capacity
  \begin{equation*}
    c-\sum_{r<t}g_{(r)}
  \end{equation*}
  left by them. Since $a_{(r)}=g_{(r)}$ for $r<t$ and $a$ is a pattern,
  $a_{(t)}$ is zero or another such item, so $a_{(t)}\le g_{(t)}$, with
  strict inequality because the sequences differ at position $t$. The
  pattern formed by $g_{(1)},\dots,g_{(t)}$ is then feasible and
  lexicographically larger than $a$, a contradiction.
  \pitem{c} First fit decreasing places the items in nonincreasing order
  into the first bin with room. Bin $k$ receives exactly the items not
  placed in bins $1,\dots,k-1$ that fit when scanned in nonincreasing
  order, which is the construction of step (b) applied to the residual
  after bins $1,\dots,k-1$. The number of residual items only decreases,
  so the condition on $s$ holds at every step of SS-$s$.
\end{proof}

\begin{proof}[Proof of Theorem~\ref{thm:smoothed}]
  \pitem{i} Consider the master with its value capped,
  \begin{equation*}
    z_U=\min\bigl\{\one^\top\lambda:\ A\lambda\ge\dv,\ \one^\top\lambda\le U,\ \lambda\ge0\bigr\}.
  \end{equation*}
  Since $U\ge\zlp$, the added constraint is satisfied by an optimal
  solution of the relaxation, so $z_U=\zlp$. The dual of this program is
  \begin{equation*}
    \max\bigl\{\dv^\top\pi-U\varsigma:\ \pi^\top a-\varsigma\le1\ \ \forall a\in\calA,\ \ \pi\ge0,\ \ \varsigma\ge0\bigr\},
  \end{equation*}
  and for fixed $\pi$ the best $\varsigma$ is $[\nu(\pi)]_+$, so the dual
  objective as a function of $\pi$ is
  \begin{equation*}
    \dv^\top\pi-U[\nu(\pi)]_+\;=\;\theta(\pi).
  \end{equation*}
  Weak duality gives $\theta\le\zlp$, and strong duality gives
  \begin{equation*}
    \sup_{\pi\ge0}\theta(\pi)\;=\;\zlp .
  \end{equation*}
  An optimal dual $\pi^*$ of the relaxation satisfies $\pi^*_i\le1$, by the
  singleton patterns, and $\nu(\pi^*)\le0$, so
  \begin{equation*}
    \theta(\pi^*)\;=\;\dv^\top\pi^*\;=\;\zlp ,
  \end{equation*}
  and the supremum is attained on $[0,1]^n$. Write $v^*=v^*(\pi)$ for the
  largest pattern value. For $\pi\ge0$ with $v^*\ge1$ the vector $\pi/v^*$
  is dual feasible, so $\dv^\top\pi\le\zlp v^*$ and
  \begin{equation*}
    \theta(\pi)\;\le\;\zlp v^*-U(v^*-1)\;=\;U-(U-\zlp)v^* .
  \end{equation*}
  For $v^*<1$ the same bound holds, since then
  \begin{equation*}
    \theta(\pi)\;\le\;\dv^\top\pi\;\le\;\zlp v^*\;<\;\zlp\;\le\;U-(U-\zlp)v^* .
  \end{equation*}
  If $U>\zlp$, a maximizer with $\nu(\pi)>0$ would have
  $\theta(\pi)<\zlp$, so the maximizers have $\nu\le0$ and
  $\dv^\top\pi=\zlp$, that is, they are the optimal duals. Since
  $v^*(\pi)\ge\max_i\pi_i$, by the singleton patterns, the bound also
  shows that $\theta$ is coercive.
  \pitem{ii} Since $e^{s[\nu]_+}=\max\bigl(1,e^{s\nu}\bigr)$,
  \begin{equation*}
    e^{s[\nu]_+}\;\le\;1+\sum_{a\in\calA}e^{s(v_a-1)}\;\le\;(1+N)\,e^{s[\nu]_+}.
  \end{equation*}
  Taking logarithms and dividing by $s$ gives
  \begin{equation*}
    [\nu]_+\;\le\;g_s\;\le\;[\nu]_++\frac{\log(1+N)}{s},
  \end{equation*}
  and multiplying by $-U$ and adding $\dv^\top\pi$ gives the bracket.
  \pitem{iii} The function $g_s$ is $1/s$ times the log-sum-exp of the
  affine functions $s(\pi^\top a-1)$, for $a\in\calA$, and $0$, hence convex
  and $C^\infty$. Its gradient and Hessian are
  \begin{equation*}
    \nabla g_s(\pi)=\E_{\tilde P_s}[a],\qquad
    \nabla^2g_s(\pi)=s\,\Cov_{\tilde P_s}(a).
  \end{equation*}
  Since $\|a\|_2\le N_{\max}$ for every pattern,
  \begin{equation*}
    \Cov_{\tilde P_s}(a)\;\preceq\;\E_{\tilde P_s}\bigl[aa^\top\bigr]\;\preceq\;N_{\max}^2\,\mathrm{Id},
  \end{equation*}
  so $\theta_s=\dv^\top\pi-Ug_s$ is concave with
  $\nabla^2\theta_s\succeq-UsN_{\max}^2\,\mathrm{Id}$, and its gradient
  $\dv-U\E_{\tilde P_s}[a]$ is Lipschitz with constant $L_s=UsN_{\max}^2$.
  \pitem{iv} The accelerated projected gradient method of Beck and
  Teboulle~\cite{beck2009}, applied to the convex $L_s$-smooth function
  $-\theta_s$ plus the indicator of $[0,1]^n$, whose projection is
  clipping, satisfies
  \begin{equation*}
    \max_{[0,1]^n}\theta_s-\theta_s(\pi^k)\;\le\;\frac{2L_s\,\|\pi^0-\pi^*_s\|^2}{(k+1)^2},
  \end{equation*}
  where $\pi^*_s$ is a maximizer of $\theta_s$ over $[0,1]^n$, and
  $\|\pi^0-\pi^*_s\|^2\le n$ on the unit box, which gives
  inequality~\eqref{eq:nesterov}. Every iterate carries the bound
  $\zlp\ge\theta(\pi^k)\ge\theta_s(\pi^k)$, by parts (i) and (ii). With
  $s$ as in condition~\eqref{eq:nesterovparams} the smoothing loss
  $U\log(1+N)/s$ equals $\epsilon/2$, and with $k$ as there the right side
  of inequality~\eqref{eq:nesterov} is at most $\epsilon/2$. For an optimal
  dual $\pi^*\in[0,1]^n$, parts (i) and (ii) then give
  \begin{equation*}
    \theta_s(\pi^k)\;\ge\;\theta_s(\pi^*)-\frac\epsilon2\;\ge\;\theta(\pi^*)-\frac{U\log(1+N)}{s}-\frac\epsilon2\;=\;\zlp-\epsilon .
  \end{equation*}
\end{proof}

\begin{proof}[Proof of Theorem~\ref{thm:sampled}]
  Fix a row $i$ and let
  \begin{equation*}
    X_i=\sum_{j=1}^ma^{(j)}_i ,
  \end{equation*}
  a sum of $m$ independent variables with values in $[0,N_{\max}]$ and
  mean
  \begin{equation*}
    \mu_i\;=\;m\,\E_{\tilde P}[a_i]\;\ge\;\frac{md_i}{U}.
  \end{equation*}
  The multiplicative Chernoff bound for independent variables with values
  in $[0,1]$~\cite{dubhashi2009}, applied to the variables $a^{(j)}_i/N_{\max}$,
  gives, for $m$ as in condition~\eqref{eq:msamples},
  \begin{equation*}
    \Prob\bigl(X_i\le(1-\xi)\mu_i\bigr)\;\le\;\exp\Bigl(-\frac{\xi^2\mu_i}{2N_{\max}}\Bigr)
    \;\le\;\exp\Bigl(-\frac{\xi^2md_{\min}}{2N_{\max}U}\Bigr)\;\le\;\frac{\zeta}{2n}.
  \end{equation*}
  On the complementary event,
  \begin{equation*}
    (A\hat\lambda^m)_i\;=\;\frac{U X_i}{m}\;>\;(1-\xi)\,U\,\E_{\tilde P}[a_i]\;\ge\;(1-\xi)\,d_i ,
  \end{equation*}
  so $\lambda'=\hat\lambda^m/(1-\xi)$ covers row $i$, and a union bound over
  the $n$ rows costs $\zeta/2$. The fraction $\hat q$ of samples that are
  nonzero patterns is a mean of $m$ Bernoulli variables with mean
  $\tilde P(\calA)$, so by Hoeffding's inequality
  \begin{equation*}
    \hat q\;<\;\tilde P(\calA)+\sqrt{\frac{\ln(2/\zeta)}{2m}}
  \end{equation*}
  with probability at least $1-\zeta/2$. On this event
  \begin{equation*}
    \one^\top\lambda'\;=\;\frac{U\hat q}{1-\xi}\;<\;\frac{1}{1-\xi}\Bigl(U\tilde P(\calA)+U\sqrt{\frac{\ln(2/\zeta)}{2m}}\Bigr),
  \end{equation*}
  which is inequality~\eqref{eq:sampledvalue}. Both events hold with
  probability at least $1-\zeta$, and the support of $\lambda'$ has at most
  $m$ columns.
\end{proof}

\begin{proof}[Proof of Lemma~\ref{lem:condexchange}]
  \pitem{i} Let a perfect packing respect $(V,J)$, let $P=\{i,x,y\}$, and
  let $B$ be the bin of $i$. It contains no item of $V$, so it is a full bin
  $F\ni i$ with $F\cap V=\emptyset$, which contains $x$ or $y$ by
  hypothesis. If it contains both, then $B=P$. Otherwise exactly one item
  $z\in\{x,y\}$ is missing from $B$, it lies in a bin $B'\ne B$, and the
  exchange in the proof of Lemma~\ref{lem:exchange} replaces the two bins
  by
  \begin{equation*}
    B\;\longmapsto\;P,\qquad B'\;\longmapsto\;(B'\setminus\{z\})\cup R',
    \qquad R'=B\setminus(P\setminus\{z\}),
  \end{equation*}
  where $R'$ weighs $w_z\le c-w_i<c/2$, since $i$ is large. The bin $P$
  contains no item of $V$. The items $R'$ contain no item of $V$, since $B$
  contains none, and no large item, since their weight is below $c/2$; so
  the new bin $B'$ contains an item of $J$ only if the old one did, and
  then the old one contained no item of $V$ and the new one contains none.
  Every other bin is unchanged, and the new packing respects $(V,J)$.
  \pitem{ii} The bin of $i$ contains no item of $V$, so its completion lies
  in $R\setminus V$ and is the one completion there, up to the interchange
  of items of equal weight.
  \pitem{iii} Let a perfect packing respect $(V,J)$, and let $B$ and $B'$ be
  the bins of $i$ and $x$. If $B=B'$, then $B=\{i,x\}$. Otherwise
  exchanging $B\setminus\{i\}$ and $x$, which both weigh $c-w_i$, replaces
  the two bins by
  \begin{equation*}
    B\;\longmapsto\;\{i,x\},\qquad B'\;\longmapsto\;L=(B\setminus\{i\})\cup(B'\setminus\{x\}),
  \end{equation*}
  both full. The bin $\{i,x\}$ contains no item of $V$. The items
  $B\setminus\{i\}$ contain no item of $V$, since the packing respects
  $(V,J)$, and no large item, since they weigh $c-w_i<c/2$; and $x$ is not
  large either. So $L$ contains an item of $J$ only if $B'$ did, and then
  $B'$ contained no item of $V$ and $L$ contains none. Every other bin is
  unchanged, and the new packing respects $(V,J)$.
\end{proof}

\end{appendices}

\bibliography{refs}


\begin{thebibliography}{65}
\ifx \bisbn   \undefined \def \bisbn  #1{ISBN #1}\fi
\ifx \binits  \undefined \def \binits#1{#1}\fi
\ifx \bauthor  \undefined \def \bauthor#1{#1}\fi
\ifx \batitle  \undefined \def \batitle#1{#1}\fi
\ifx \bjtitle  \undefined \def \bjtitle#1{#1}\fi
\ifx \bvolume  \undefined \def \bvolume#1{\textbf{#1}}\fi
\ifx \byear  \undefined \def \byear#1{#1}\fi
\ifx \bissue  \undefined \def \bissue#1{#1}\fi
\ifx \bfpage  \undefined \def \bfpage#1{#1}\fi
\ifx \blpage  \undefined \def \blpage #1{#1}\fi
\ifx \burl  \undefined \def \burl#1{\textsf{#1}}\fi
\ifx \doiurl  \undefined \def \doiurl#1{\url{https://doi.org/#1}}\fi
\ifx \betal  \undefined \def \betal{\textit{et al.}}\fi
\ifx \binstitute  \undefined \def \binstitute#1{#1}\fi
\ifx \binstitutionaled  \undefined \def \binstitutionaled#1{#1}\fi
\ifx \bctitle  \undefined \def \bctitle#1{#1}\fi
\ifx \beditor  \undefined \def \beditor#1{#1}\fi
\ifx \bpublisher  \undefined \def \bpublisher#1{#1}\fi
\ifx \bbtitle  \undefined \def \bbtitle#1{#1}\fi
\ifx \bedition  \undefined \def \bedition#1{#1}\fi
\ifx \bseriesno  \undefined \def \bseriesno#1{#1}\fi
\ifx \blocation  \undefined \def \blocation#1{#1}\fi
\ifx \bsertitle  \undefined \def \bsertitle#1{#1}\fi
\ifx \bsnm \undefined \def \bsnm#1{#1}\fi
\ifx \bsuffix \undefined \def \bsuffix#1{#1}\fi
\ifx \bparticle \undefined \def \bparticle#1{#1}\fi
\ifx \barticle \undefined \def \barticle#1{#1}\fi
\bibcommenthead
\ifx \bconfdate \undefined \def \bconfdate #1{#1}\fi
\ifx \botherref \undefined \def \botherref #1{#1}\fi
\ifx \url \undefined \def \url#1{\textsf{#1}}\fi
\ifx \bchapter \undefined \def \bchapter#1{#1}\fi
\ifx \bbook \undefined \def \bbook#1{#1}\fi
\ifx \bcomment \undefined \def \bcomment#1{#1}\fi
\ifx \oauthor \undefined \def \oauthor#1{#1}\fi
\ifx \citeauthoryear \undefined \def \citeauthoryear#1{#1}\fi
\ifx \endbibitem  \undefined \def \endbibitem {}\fi
\ifx \bconflocation  \undefined \def \bconflocation#1{#1}\fi
\ifx \arxivurl  \undefined \def \arxivurl#1{\textsf{#1}}\fi
\csname PreBibitemsHook\endcsname

\bibitem[\protect\citeauthoryear{Gilmore and Gomory}{1961}]{gilmore1961}
\begin{barticle}
\bauthor{\bsnm{Gilmore}, \binits{P.C.}},
\bauthor{\bsnm{Gomory}, \binits{R.E.}}:
\batitle{A linear programming approach to the cutting-stock problem}.
\bjtitle{Operations Research}
\bvolume{9}(\bissue{6}),
\bfpage{849}--\blpage{859}
(\byear{1961})
\end{barticle}
\endbibitem

\bibitem[\protect\citeauthoryear{Uchoa and Sadykov}{2026}]{uchoa2026}
\begin{barticle}
\bauthor{\bsnm{Uchoa}, \binits{E.}},
\bauthor{\bsnm{Sadykov}, \binits{R.}}:
\batitle{Kantorovich and {Z}algaller (1951): the 0-th column generation
  algorithm}.
\bjtitle{Mathematical Programming}
\bvolume{217},
\bfpage{439}--\blpage{456}
(\byear{2026})
\end{barticle}
\endbibitem

\bibitem[\protect\citeauthoryear{Scheithauer and Terno}{1995}]{scheithauer1995}
\begin{barticle}
\bauthor{\bsnm{Scheithauer}, \binits{G.}},
\bauthor{\bsnm{Terno}, \binits{J.}}:
\batitle{The modified integer round-up property of the one-dimensional cutting
  stock problem}.
\bjtitle{European Journal of Operational Research}
\bvolume{84}(\bissue{3}),
\bfpage{562}--\blpage{571}
(\byear{1995})
\end{barticle}
\endbibitem

\bibitem[\protect\citeauthoryear{Kartak et~al.}{2015}]{kartak2015}
\begin{barticle}
\bauthor{\bsnm{Kartak}, \binits{V.M.}},
\bauthor{\bsnm{Ripatti}, \binits{A.V.}},
\bauthor{\bsnm{Scheithauer}, \binits{G.}},
\bauthor{\bsnm{Kurz}, \binits{S.}}:
\batitle{Minimal proper non-{IRUP} instances of the one-dimensional cutting
  stock problem}.
\bjtitle{Discrete Applied Mathematics}
\bvolume{187},
\bfpage{120}--\blpage{129}
(\byear{2015})
\end{barticle}
\endbibitem

\bibitem[\protect\citeauthoryear{Goemans and Rothvoss}{2020}]{goemans2020}
\begin{barticle}
\bauthor{\bsnm{Goemans}, \binits{M.X.}},
\bauthor{\bsnm{Rothvoss}, \binits{T.}}:
\batitle{Polynomiality for bin packing with a constant number of item types}.
\bjtitle{Journal of the ACM}
\bvolume{67}(\bissue{6}),
\bfpage{1}--\blpage{21}
(\byear{2020}).
\bcomment{Article 38}
\end{barticle}
\endbibitem

\bibitem[\protect\citeauthoryear{Vance et~al.}{1994}]{vance1994}
\begin{barticle}
\bauthor{\bsnm{Vance}, \binits{P.H.}},
\bauthor{\bsnm{Barnhart}, \binits{C.}},
\bauthor{\bsnm{Johnson}, \binits{E.L.}},
\bauthor{\bsnm{Nemhauser}, \binits{G.L.}}:
\batitle{Solving binary cutting stock problems by column generation and
  branch-and-bound}.
\bjtitle{Computational Optimization and Applications}
\bvolume{3}(\bissue{2}),
\bfpage{111}--\blpage{130}
(\byear{1994})
\end{barticle}
\endbibitem

\bibitem[\protect\citeauthoryear{Vanderbeck}{1999}]{vanderbeck1999}
\begin{barticle}
\bauthor{\bsnm{Vanderbeck}, \binits{F.}}:
\batitle{Computational study of a column generation algorithm for bin packing
  and cutting stock problems}.
\bjtitle{Mathematical Programming}
\bvolume{86}(\bissue{3}),
\bfpage{565}--\blpage{594}
(\byear{1999})
\end{barticle}
\endbibitem

\bibitem[\protect\citeauthoryear{Barnhart et~al.}{1998}]{barnhart1998}
\begin{barticle}
\bauthor{\bsnm{Barnhart}, \binits{C.}},
\bauthor{\bsnm{Johnson}, \binits{E.L.}},
\bauthor{\bsnm{Nemhauser}, \binits{G.L.}},
\bauthor{\bsnm{Savelsbergh}, \binits{M.W.P.}},
\bauthor{\bsnm{Vance}, \binits{P.H.}}:
\batitle{Branch-and-price: column generation for solving huge integer
  programs}.
\bjtitle{Operations Research}
\bvolume{46}(\bissue{3}),
\bfpage{316}--\blpage{329}
(\byear{1998})
\end{barticle}
\endbibitem

\bibitem[\protect\citeauthoryear{Belov and Scheithauer}{2006}]{belov2006}
\begin{barticle}
\bauthor{\bsnm{Belov}, \binits{G.}},
\bauthor{\bsnm{Scheithauer}, \binits{G.}}:
\batitle{A branch-and-cut-and-price algorithm for one-dimensional stock cutting
  and two-dimensional two-stage cutting}.
\bjtitle{European Journal of Operational Research}
\bvolume{171}(\bissue{1}),
\bfpage{85}--\blpage{106}
(\byear{2006})
\end{barticle}
\endbibitem

\bibitem[\protect\citeauthoryear{{Val{\'e}rio de Carvalho}}{1999}]{valerio1999}
\begin{barticle}
\bauthor{\bsnm{{Val{\'e}rio de Carvalho}}, \binits{J.M.}}:
\batitle{Exact solution of bin-packing problems using column generation and
  branch-and-bound}.
\bjtitle{Annals of Operations Research}
\bvolume{86},
\bfpage{629}--\blpage{659}
(\byear{1999})
\end{barticle}
\endbibitem

\bibitem[\protect\citeauthoryear{Brand{\~a}o and Pedroso}{2016}]{brandao2016}
\begin{barticle}
\bauthor{\bsnm{Brand{\~a}o}, \binits{F.}},
\bauthor{\bsnm{Pedroso}, \binits{J.P.}}:
\batitle{Bin packing and related problems: general arc-flow formulation with
  graph compression}.
\bjtitle{Computers \& Operations Research}
\bvolume{69},
\bfpage{56}--\blpage{67}
(\byear{2016})
\end{barticle}
\endbibitem

\bibitem[\protect\citeauthoryear{{de Lima} et~al.}{2022}]{delima2022}
\begin{barticle}
\bauthor{\bsnm{{de Lima}}, \binits{V.L.}},
\bauthor{\bsnm{Alves}, \binits{C.}},
\bauthor{\bsnm{Clautiaux}, \binits{F.}},
\bauthor{\bsnm{Iori}, \binits{M.}},
\bauthor{\bsnm{{Val{\'e}rio de Carvalho}}, \binits{J.M.}}:
\batitle{Arc flow formulations based on dynamic programming: theoretical
  foundations and applications}.
\bjtitle{European Journal of Operational Research}
\bvolume{296}(\bissue{1}),
\bfpage{3}--\blpage{21}
(\byear{2022})
\end{barticle}
\endbibitem

\bibitem[\protect\citeauthoryear{Delorme and Iori}{2020}]{delorme2020}
\begin{barticle}
\bauthor{\bsnm{Delorme}, \binits{M.}},
\bauthor{\bsnm{Iori}, \binits{M.}}:
\batitle{Enhanced pseudo-polynomial formulations for bin packing and cutting
  stock problems}.
\bjtitle{INFORMS Journal on Computing}
\bvolume{32}(\bissue{1}),
\bfpage{101}--\blpage{119}
(\byear{2020})
\end{barticle}
\endbibitem

\bibitem[\protect\citeauthoryear{Delorme and Martinovic}{2026}]{delorme2026}
\begin{barticle}
\bauthor{\bsnm{Delorme}, \binits{M.}},
\bauthor{\bsnm{Martinovic}, \binits{J.}}:
\batitle{Extending the reflect flow formulation to variable-sized
  one-dimensional cutting and skiving stock problems}.
\bjtitle{European Journal of Operational Research}
\bvolume{334},
\bfpage{51}--\blpage{69}
(\byear{2026})
\end{barticle}
\endbibitem

\bibitem[\protect\citeauthoryear{Delorme et~al.}{2016}]{delorme2016}
\begin{barticle}
\bauthor{\bsnm{Delorme}, \binits{M.}},
\bauthor{\bsnm{Iori}, \binits{M.}},
\bauthor{\bsnm{Martello}, \binits{S.}}:
\batitle{Bin packing and cutting stock problems: Mathematical models and exact
  algorithms}.
\bjtitle{European Journal of Operational Research}
\bvolume{255}(\bissue{1}),
\bfpage{1}--\blpage{20}
(\byear{2016})
\end{barticle}
\endbibitem

\bibitem[\protect\citeauthoryear{Ataei and Wang}{2022}]{ataei2022}
\begin{barticle}
\bauthor{\bsnm{Ataei}, \binits{M.}},
\bauthor{\bsnm{Wang}, \binits{X.}}:
\batitle{Theory of {L}ehmer transform and its applications in identifying the
  electroencephalographic signature of major depressive disorder}.
\bjtitle{Scientific Reports}
\bvolume{12}(\bissue{1}),
\bfpage{3663}
(\byear{2022})
\end{barticle}
\endbibitem

\bibitem[\protect\citeauthoryear{Ataei}{2026}]{ataei2026}
\begin{botherref}
\oauthor{\bsnm{Ataei}, \binits{M.}}:
{L}ehmer transform.
arXiv:2609.05790
(2026)
\end{botherref}
\endbibitem

\bibitem[\protect\citeauthoryear{Ataei and
  Vinogradov}{2026}]{ataei_vinogradov2026}
\begin{botherref}
\oauthor{\bsnm{Ataei}, \binits{M.}},
\oauthor{\bsnm{Vinogradov}, \binits{V.V.}}:
{L}ehmer--{L}ambert distribution.
arXiv:2609.12150
(2026)
\end{botherref}
\endbibitem

\bibitem[\protect\citeauthoryear{Caprara et~al.}{2015}]{caprara2015}
\begin{barticle}
\bauthor{\bsnm{Caprara}, \binits{A.}},
\bauthor{\bsnm{Dell'Amico}, \binits{M.}},
\bauthor{\bsnm{D{\'\i}az-D{\'\i}az}, \binits{J.C.}},
\bauthor{\bsnm{Iori}, \binits{M.}},
\bauthor{\bsnm{Rizzi}, \binits{R.}}:
\batitle{Friendly bin packing instances without integer round-up property}.
\bjtitle{Mathematical Programming}
\bvolume{150}(\bissue{1}),
\bfpage{5}--\blpage{17}
(\byear{2015})
\end{barticle}
\endbibitem

\bibitem[\protect\citeauthoryear{Wei et~al.}{2020}]{wei2020}
\begin{barticle}
\bauthor{\bsnm{Wei}, \binits{L.}},
\bauthor{\bsnm{Luo}, \binits{Z.}},
\bauthor{\bsnm{Baldacci}, \binits{R.}},
\bauthor{\bsnm{Lim}, \binits{A.}}:
\batitle{A new branch-and-price-and-cut algorithm for one-dimensional
  bin-packing problems}.
\bjtitle{INFORMS Journal on Computing}
\bvolume{32}(\bissue{2}),
\bfpage{428}--\blpage{443}
(\byear{2020})
\end{barticle}
\endbibitem

\bibitem[\protect\citeauthoryear{Pessoa et~al.}{2020}]{pessoa2020}
\begin{barticle}
\bauthor{\bsnm{Pessoa}, \binits{A.}},
\bauthor{\bsnm{Sadykov}, \binits{R.}},
\bauthor{\bsnm{Uchoa}, \binits{E.}},
\bauthor{\bsnm{Vanderbeck}, \binits{F.}}:
\batitle{A generic exact solver for vehicle routing and related problems}.
\bjtitle{Mathematical Programming}
\bvolume{183},
\bfpage{483}--\blpage{523}
(\byear{2020})
\end{barticle}
\endbibitem

\bibitem[\protect\citeauthoryear{{de Lima} et~al.}{2023}]{delima2023}
\begin{barticle}
\bauthor{\bsnm{{de Lima}}, \binits{V.L.}},
\bauthor{\bsnm{Iori}, \binits{M.}},
\bauthor{\bsnm{Miyazawa}, \binits{F.K.}}:
\batitle{Exact solution of network flow models with strong relaxations}.
\bjtitle{Mathematical Programming}
\bvolume{197},
\bfpage{813}--\blpage{846}
(\byear{2023})
\end{barticle}
\endbibitem

\bibitem[\protect\citeauthoryear{{da Silva} and Schouery}{2025}]{dasilva2025}
\begin{barticle}
\bauthor{\bsnm{{da Silva}}, \binits{R.F.F.}},
\bauthor{\bsnm{Schouery}, \binits{R.C.S.}}:
\batitle{Solving cutting stock problems via an extended {R}yan--{F}oster
  branching scheme and fast column generation}.
\bjtitle{INFORMS Journal on Computing}
(\byear{2025})
\doiurl{10.1287/ijoc.2023.0399}
\end{barticle}
\endbibitem

\bibitem[\protect\citeauthoryear{{da Silva} et~al.}{2026}]{dasilva2026}
\begin{botherref}
\oauthor{\bsnm{{da Silva}}, \binits{R.F.F.}},
\oauthor{\bsnm{{de Lima}}, \binits{V.L.}},
\oauthor{\bsnm{Schouery}, \binits{R.C.S.}},
\oauthor{\bsnm{C\^ot\'e}, \binits{J.-F.}},
\oauthor{\bsnm{Iori}, \binits{M.}}:
Polynomial and pseudopolynomial algorithms for two classes of bin packing
  instances.
arXiv:2604.05152
(2026)
\end{botherref}
\endbibitem

\bibitem[\protect\citeauthoryear{Williams et~al.}{2003}]{williams2003}
\begin{bchapter}
\bauthor{\bsnm{Williams}, \binits{R.}},
\bauthor{\bsnm{Gomes}, \binits{C.P.}},
\bauthor{\bsnm{Selman}, \binits{B.}}:
\bctitle{Backdoors to typical case complexity}.
In: \bbtitle{Proceedings of the 18th International Joint Conference on
  Artificial Intelligence ({IJCAI})},
pp. \bfpage{1173}--\blpage{1178}
(\byear{2003})
\end{bchapter}
\endbibitem

\bibitem[\protect\citeauthoryear{Fischetti and Monaci}{2013}]{fischetti2013}
\begin{barticle}
\bauthor{\bsnm{Fischetti}, \binits{M.}},
\bauthor{\bsnm{Monaci}, \binits{M.}}:
\batitle{Backdoor branching}.
\bjtitle{INFORMS Journal on Computing}
\bvolume{25}(\bissue{4}),
\bfpage{693}--\blpage{700}
(\byear{2013})
\end{barticle}
\endbibitem

\bibitem[\protect\citeauthoryear{Lehmer}{1971}]{lehmer1971}
\begin{barticle}
\bauthor{\bsnm{Lehmer}, \binits{D.H.}}:
\batitle{On the compounding of certain means}.
\bjtitle{Journal of Mathematical Analysis and Applications}
\bvolume{36}(\bissue{1}),
\bfpage{183}--\blpage{200}
(\byear{1971})
\end{barticle}
\endbibitem

\bibitem[\protect\citeauthoryear{Bullen}{2003}]{bullen2003}
\begin{bbook}
\bauthor{\bsnm{Bullen}, \binits{P.S.}}:
\bbtitle{Handbook of Means and Their Inequalities}.
\bpublisher{Kluwer Academic Publishers},
\blocation{Dordrecht}
(\byear{2003})
\end{bbook}
\endbibitem

\bibitem[\protect\citeauthoryear{Clautiaux et~al.}{2010}]{clautiaux2010}
\begin{barticle}
\bauthor{\bsnm{Clautiaux}, \binits{F.}},
\bauthor{\bsnm{Alves}, \binits{C.}},
\bauthor{\bsnm{{Val{\'e}rio de Carvalho}}, \binits{J.}}:
\batitle{A survey of dual-feasible and superadditive functions}.
\bjtitle{Annals of Operations Research}
\bvolume{179},
\bfpage{317}--\blpage{342}
(\byear{2010})
\end{barticle}
\endbibitem

\bibitem[\protect\citeauthoryear{L{\"u}bbecke and
  Desrosiers}{2005}]{lubbecke2005}
\begin{barticle}
\bauthor{\bsnm{L{\"u}bbecke}, \binits{M.E.}},
\bauthor{\bsnm{Desrosiers}, \binits{J.}}:
\batitle{Selected topics in column generation}.
\bjtitle{Operations Research}
\bvolume{53}(\bissue{6}),
\bfpage{1007}--\blpage{1023}
(\byear{2005})
\end{barticle}
\endbibitem

\bibitem[\protect\citeauthoryear{Held et~al.}{2012}]{heldcook2012}
\begin{barticle}
\bauthor{\bsnm{Held}, \binits{S.}},
\bauthor{\bsnm{Cook}, \binits{W.}},
\bauthor{\bsnm{Sewell}, \binits{E.C.}}:
\batitle{Maximum-weight stable sets and safe lower bounds for graph coloring}.
\bjtitle{Mathematical Programming Computation}
\bvolume{4}(\bissue{4}),
\bfpage{363}--\blpage{381}
(\byear{2012})
\end{barticle}
\endbibitem

\bibitem[\protect\citeauthoryear{Baldacci et~al.}{2024}]{baldacci2024}
\begin{barticle}
\bauthor{\bsnm{Baldacci}, \binits{R.}},
\bauthor{\bsnm{Coniglio}, \binits{S.}},
\bauthor{\bsnm{Cordeau}, \binits{J.-F.}},
\bauthor{\bsnm{Furini}, \binits{F.}}:
\batitle{A numerically exact algorithm for the bin-packing problem}.
\bjtitle{INFORMS Journal on Computing}
\bvolume{36}(\bissue{1}),
\bfpage{141}--\blpage{162}
(\byear{2024})
\end{barticle}
\endbibitem

\bibitem[\protect\citeauthoryear{Cominetti and {San
  Mart{\'i}n}}{1994}]{cominetti1994}
\begin{barticle}
\bauthor{\bsnm{Cominetti}, \binits{R.}},
\bauthor{\bsnm{{San Mart{\'i}n}}, \binits{J.}}:
\batitle{Asymptotic analysis of the exponential penalty trajectory in linear
  programming}.
\bjtitle{Mathematical Programming}
\bvolume{67}(\bissue{1--3}),
\bfpage{169}--\blpage{187}
(\byear{1994})
\end{barticle}
\endbibitem

\bibitem[\protect\citeauthoryear{Nesterov}{2005}]{nesterov2005}
\begin{barticle}
\bauthor{\bsnm{Nesterov}, \binits{Y.}}:
\batitle{Smooth minimization of non-smooth functions}.
\bjtitle{Mathematical Programming}
\bvolume{103}(\bissue{1}),
\bfpage{127}--\blpage{152}
(\byear{2005})
\end{barticle}
\endbibitem

\bibitem[\protect\citeauthoryear{Bienstock and Iyengar}{2006}]{bienstock2006}
\begin{barticle}
\bauthor{\bsnm{Bienstock}, \binits{D.}},
\bauthor{\bsnm{Iyengar}, \binits{G.}}:
\batitle{Approximating fractional packings and coverings in {$O(1/\epsilon)$}
  iterations}.
\bjtitle{SIAM Journal on Computing}
\bvolume{35}(\bissue{4}),
\bfpage{825}--\blpage{854}
(\byear{2006})
\end{barticle}
\endbibitem

\bibitem[\protect\citeauthoryear{Martello and Toth}{1990}]{martello1990}
\begin{bbook}
\bauthor{\bsnm{Martello}, \binits{S.}},
\bauthor{\bsnm{Toth}, \binits{P.}}:
\bbtitle{Knapsack Problems: Algorithms and Computer Implementations}.
\bpublisher{John Wiley \& Sons},
\blocation{Chichester}
(\byear{1990})
\end{bbook}
\endbibitem

\bibitem[\protect\citeauthoryear{Gupta and Ho}{1999}]{gupta1999}
\begin{barticle}
\bauthor{\bsnm{Gupta}, \binits{J.N.D.}},
\bauthor{\bsnm{Ho}, \binits{J.C.}}:
\batitle{A new heuristic algorithm for the one-dimensional bin-packing
  problem}.
\bjtitle{Production Planning \& Control}
\bvolume{10}(\bissue{6}),
\bfpage{598}--\blpage{603}
(\byear{1999})
\end{barticle}
\endbibitem

\bibitem[\protect\citeauthoryear{Fleszar and Hindi}{2002}]{fleszar2002}
\begin{barticle}
\bauthor{\bsnm{Fleszar}, \binits{K.}},
\bauthor{\bsnm{Hindi}, \binits{K.S.}}:
\batitle{New heuristics for one-dimensional bin-packing}.
\bjtitle{Computers \& Operations Research}
\bvolume{29}(\bissue{7}),
\bfpage{821}--\blpage{839}
(\byear{2002})
\end{barticle}
\endbibitem

\bibitem[\protect\citeauthoryear{Quiroz-Castellanos et~al.}{2015}]{quiroz2015}
\begin{barticle}
\bauthor{\bsnm{Quiroz-Castellanos}, \binits{M.}},
\bauthor{\bsnm{Cruz-Reyes}, \binits{L.}},
\bauthor{\bsnm{Torres-Jim{\'e}nez}, \binits{J.}},
\bauthor{\bsnm{G{\'o}mez~Santill{\'a}n}, \binits{C.}},
\bauthor{\bsnm{Fraire~Huacuja}, \binits{H.J.}},
\bauthor{\bsnm{Alvim}, \binits{A.C.F.}}:
\batitle{A grouping genetic algorithm with controlled gene transmission for the
  bin packing problem}.
\bjtitle{Computers \& Operations Research}
\bvolume{55},
\bfpage{52}--\blpage{64}
(\byear{2015})
\end{barticle}
\endbibitem

\bibitem[\protect\citeauthoryear{Caprara and Pferschy}{2004}]{caprara2004}
\begin{barticle}
\bauthor{\bsnm{Caprara}, \binits{A.}},
\bauthor{\bsnm{Pferschy}, \binits{U.}}:
\batitle{Worst-case analysis of the subset sum algorithm for bin packing}.
\bjtitle{Operations Research Letters}
\bvolume{32}(\bissue{2}),
\bfpage{159}--\blpage{166}
(\byear{2004})
\end{barticle}
\endbibitem

\bibitem[\protect\citeauthoryear{Caprara and Pferschy}{2005}]{caprara2005}
\begin{barticle}
\bauthor{\bsnm{Caprara}, \binits{A.}},
\bauthor{\bsnm{Pferschy}, \binits{U.}}:
\batitle{Modified subset sum heuristics for bin packing}.
\bjtitle{Information Processing Letters}
\bvolume{96}(\bissue{1}),
\bfpage{18}--\blpage{23}
(\byear{2005})
\end{barticle}
\endbibitem

\bibitem[\protect\citeauthoryear{Farley}{1990}]{farley1990}
\begin{barticle}
\bauthor{\bsnm{Farley}, \binits{A.A.}}:
\batitle{A note on bounding a class of linear programming problems, including
  cutting stock problems}.
\bjtitle{Operations Research}
\bvolume{38}(\bissue{5}),
\bfpage{922}--\blpage{923}
(\byear{1990})
\end{barticle}
\endbibitem

\bibitem[\protect\citeauthoryear{Garg and K{\"o}nemann}{2007}]{garg2007}
\begin{barticle}
\bauthor{\bsnm{Garg}, \binits{N.}},
\bauthor{\bsnm{K{\"o}nemann}, \binits{J.}}:
\batitle{Faster and simpler algorithms for multicommodity flow and other
  fractional packing problems}.
\bjtitle{SIAM Journal on Computing}
\bvolume{37}(\bissue{2}),
\bfpage{630}--\blpage{652}
(\byear{2007})
\end{barticle}
\endbibitem

\bibitem[\protect\citeauthoryear{Dilkina et~al.}{2009}]{dilkina2009}
\begin{bchapter}
\bauthor{\bsnm{Dilkina}, \binits{B.}},
\bauthor{\bsnm{Gomes}, \binits{C.P.}},
\bauthor{\bsnm{Malitsky}, \binits{Y.}},
\bauthor{\bsnm{Sabharwal}, \binits{A.}},
\bauthor{\bsnm{Sellmann}, \binits{M.}}:
\bctitle{Backdoors to combinatorial optimization: feasibility and optimality}.
In: \bbtitle{Integration of {AI} and {OR} Techniques in Constraint Programming
  for Combinatorial Optimization Problems ({CPAIOR} 2009)}.
\bsertitle{Lecture Notes in Computer Science},
vol. \bseriesno{5547},
pp. \bfpage{56}--\blpage{70}.
\bpublisher{Springer}, \blocation{???}
(\byear{2009})
\end{bchapter}
\endbibitem

\bibitem[\protect\citeauthoryear{Khalil et~al.}{2022}]{khalil2022}
\begin{bchapter}
\bauthor{\bsnm{Khalil}, \binits{E.B.}},
\bauthor{\bsnm{Vaezipoor}, \binits{P.}},
\bauthor{\bsnm{Dilkina}, \binits{B.}}:
\bctitle{Finding backdoors to integer programs: a {M}onte {C}arlo tree search
  framework}.
In: \bbtitle{Proceedings of the {AAAI} Conference on Artificial Intelligence},
vol. \bseriesno{36}
(\byear{2022})
\end{bchapter}
\endbibitem

\bibitem[\protect\citeauthoryear{Higham}{2002}]{higham2002}
\begin{bbook}
\bauthor{\bsnm{Higham}, \binits{N.J.}}:
\bbtitle{Accuracy and Stability of Numerical Algorithms},
\bedition{2}nd edn.
\bpublisher{SIAM},
\blocation{Philadelphia}
(\byear{2002})
\end{bbook}
\endbibitem

\bibitem[\protect\citeauthoryear{Ryan and Foster}{1981}]{ryan1981}
\begin{bchapter}
\bauthor{\bsnm{Ryan}, \binits{D.M.}},
\bauthor{\bsnm{Foster}, \binits{B.A.}}:
\bctitle{An integer programming approach to scheduling}.
In: \beditor{\bsnm{Wren}, \binits{A.}} (ed.)
\bbtitle{Computer Scheduling of Public Transport: Urban Passenger Vehicle and
  Crew Scheduling},
pp. \bfpage{269}--\blpage{280}.
\bpublisher{North-Holland},
\blocation{Amsterdam}
(\byear{1981})
\end{bchapter}
\endbibitem

\bibitem[\protect\citeauthoryear{Vanderbeck}{2011}]{vanderbeck2011}
\begin{barticle}
\bauthor{\bsnm{Vanderbeck}, \binits{F.}}:
\batitle{Branching in branch-and-price: a generic scheme}.
\bjtitle{Mathematical Programming}
\bvolume{130}(\bissue{2}),
\bfpage{249}--\blpage{294}
(\byear{2011})
\end{barticle}
\endbibitem

\bibitem[\protect\citeauthoryear{Jepsen et~al.}{2008}]{jepsen2008}
\begin{barticle}
\bauthor{\bsnm{Jepsen}, \binits{M.}},
\bauthor{\bsnm{Petersen}, \binits{B.}},
\bauthor{\bsnm{Spoorendonk}, \binits{S.}},
\bauthor{\bsnm{Pisinger}, \binits{D.}}:
\batitle{Subset-row inequalities applied to the vehicle-routing problem with
  time windows}.
\bjtitle{Operations Research}
\bvolume{56}(\bissue{2}),
\bfpage{497}--\blpage{511}
(\byear{2008})
\end{barticle}
\endbibitem

\bibitem[\protect\citeauthoryear{Martello et~al.}{1999}]{martello1999}
\begin{barticle}
\bauthor{\bsnm{Martello}, \binits{S.}},
\bauthor{\bsnm{Pisinger}, \binits{D.}},
\bauthor{\bsnm{Toth}, \binits{P.}}:
\batitle{Dynamic programming and strong bounds for the 0-1 knapsack problem}.
\bjtitle{Management Science}
\bvolume{45}(\bissue{3}),
\bfpage{414}--\blpage{424}
(\byear{1999})
\end{barticle}
\endbibitem

\bibitem[\protect\citeauthoryear{Wentges}{1997}]{wentges1997}
\begin{barticle}
\bauthor{\bsnm{Wentges}, \binits{P.}}:
\batitle{Weighted {D}antzig--{W}olfe decomposition for linear mixed-integer
  programming}.
\bjtitle{International Transactions in Operational Research}
\bvolume{4}(\bissue{2}),
\bfpage{151}--\blpage{162}
(\byear{1997})
\end{barticle}
\endbibitem

\bibitem[\protect\citeauthoryear{Pessoa et~al.}{2018}]{pessoa2018}
\begin{barticle}
\bauthor{\bsnm{Pessoa}, \binits{A.}},
\bauthor{\bsnm{Sadykov}, \binits{R.}},
\bauthor{\bsnm{Uchoa}, \binits{E.}},
\bauthor{\bsnm{Vanderbeck}, \binits{F.}}:
\batitle{Automation and combination of linear-programming based stabilization
  techniques in column generation}.
\bjtitle{INFORMS Journal on Computing}
\bvolume{30}(\bissue{2}),
\bfpage{339}--\blpage{360}
(\byear{2018})
\end{barticle}
\endbibitem

\bibitem[\protect\citeauthoryear{{Val{\'e}rio de Carvalho}}{2005}]{valerio2005}
\begin{barticle}
\bauthor{\bsnm{{Val{\'e}rio de Carvalho}}, \binits{J.M.}}:
\batitle{Using extra dual cuts to accelerate convergence in column generation}.
\bjtitle{INFORMS Journal on Computing}
\bvolume{17}(\bissue{2}),
\bfpage{175}--\blpage{182}
(\byear{2005})
\end{barticle}
\endbibitem

\bibitem[\protect\citeauthoryear{{Ben Amor} et~al.}{2006}]{benamor2006}
\begin{barticle}
\bauthor{\bsnm{{Ben Amor}}, \binits{H.}},
\bauthor{\bsnm{Desrosiers}, \binits{J.}},
\bauthor{\bsnm{{Val{\'e}rio de Carvalho}}, \binits{J.M.}}:
\batitle{Dual-optimal inequalities for stabilized column generation}.
\bjtitle{Operations Research}
\bvolume{54}(\bissue{3}),
\bfpage{454}--\blpage{463}
(\byear{2006})
\end{barticle}
\endbibitem

\bibitem[\protect\citeauthoryear{Delorme et~al.}{2018}]{bpplib}
\begin{barticle}
\bauthor{\bsnm{Delorme}, \binits{M.}},
\bauthor{\bsnm{Iori}, \binits{M.}},
\bauthor{\bsnm{Martello}, \binits{S.}}:
\batitle{{BPPLIB}: a library for bin packing and cutting stock problems}.
\bjtitle{Optimization Letters}
\bvolume{12}(\bissue{2}),
\bfpage{235}--\blpage{250}
(\byear{2018})
\end{barticle}
\endbibitem

\bibitem[\protect\citeauthoryear{Gschwind and Irnich}{2016}]{gschwind2016}
\begin{barticle}
\bauthor{\bsnm{Gschwind}, \binits{T.}},
\bauthor{\bsnm{Irnich}, \binits{S.}}:
\batitle{Dual inequalities for stabilized column generation revisited}.
\bjtitle{INFORMS Journal on Computing}
\bvolume{28}(\bissue{1}),
\bfpage{175}--\blpage{194}
(\byear{2016})
\end{barticle}
\endbibitem

\bibitem[\protect\citeauthoryear{Falkenauer}{1996}]{falkenauer1996}
\begin{barticle}
\bauthor{\bsnm{Falkenauer}, \binits{E.}}:
\batitle{A hybrid grouping genetic algorithm for bin packing}.
\bjtitle{Journal of Heuristics}
\bvolume{2}(\bissue{1}),
\bfpage{5}--\blpage{30}
(\byear{1996})
\end{barticle}
\endbibitem

\bibitem[\protect\citeauthoryear{Scholl et~al.}{1997}]{scholl1997}
\begin{barticle}
\bauthor{\bsnm{Scholl}, \binits{A.}},
\bauthor{\bsnm{Klein}, \binits{R.}},
\bauthor{\bsnm{J{\"u}rgens}, \binits{C.}}:
\batitle{{BISON}: a fast hybrid procedure for exactly solving the
  one-dimensional bin packing problem}.
\bjtitle{Computers \& Operations Research}
\bvolume{24}(\bissue{7}),
\bfpage{627}--\blpage{645}
(\byear{1997})
\end{barticle}
\endbibitem

\bibitem[\protect\citeauthoryear{Schwerin and W{\"a}scher}{1997}]{schwerin1997}
\begin{barticle}
\bauthor{\bsnm{Schwerin}, \binits{P.}},
\bauthor{\bsnm{W{\"a}scher}, \binits{G.}}:
\batitle{The bin-packing problem: a problem generator and some numerical
  experiments with {FFD} packing and {MTP}}.
\bjtitle{International Transactions in Operational Research}
\bvolume{4}(\bissue{5--6}),
\bfpage{377}--\blpage{389}
(\byear{1997})
\end{barticle}
\endbibitem

\bibitem[\protect\citeauthoryear{W{\"a}scher and Gau}{1996}]{wascher1996}
\begin{barticle}
\bauthor{\bsnm{W{\"a}scher}, \binits{G.}},
\bauthor{\bsnm{Gau}, \binits{T.}}:
\batitle{Heuristics for the integer one-dimensional cutting stock problem: a
  computational study}.
\bjtitle{OR Spektrum}
\bvolume{18}(\bissue{3}),
\bfpage{131}--\blpage{144}
(\byear{1996})
\end{barticle}
\endbibitem

\bibitem[\protect\citeauthoryear{Schoenfield}{2002}]{schoenfield2002}
\begin{botherref}
\oauthor{\bsnm{Schoenfield}, \binits{J.E.}}:
Fast, exact solution of open bin packing problems without linear programming.
Draft, US Army Space and Missile Defense Command, Huntsville, Alabama
(2002)
\end{botherref}
\endbibitem

\bibitem[\protect\citeauthoryear{Pisinger}{2000}]{pisinger2000}
\begin{barticle}
\bauthor{\bsnm{Pisinger}, \binits{D.}}:
\batitle{A minimal algorithm for the bounded knapsack problem}.
\bjtitle{INFORMS Journal on Computing}
\bvolume{12}(\bissue{1}),
\bfpage{75}--\blpage{82}
(\byear{2000})
\end{barticle}
\endbibitem

\bibitem[\protect\citeauthoryear{Huangfu and Hall}{2018}]{huangfu2018}
\begin{barticle}
\bauthor{\bsnm{Huangfu}, \binits{Q.}},
\bauthor{\bsnm{Hall}, \binits{J.A.J.}}:
\batitle{Parallelizing the dual revised simplex method}.
\bjtitle{Mathematical Programming Computation}
\bvolume{10}(\bissue{1}),
\bfpage{119}--\blpage{142}
(\byear{2018})
\end{barticle}
\endbibitem

\bibitem[\protect\citeauthoryear{Beck and Teboulle}{2009}]{beck2009}
\begin{barticle}
\bauthor{\bsnm{Beck}, \binits{A.}},
\bauthor{\bsnm{Teboulle}, \binits{M.}}:
\batitle{A fast iterative shrinkage-thresholding algorithm for linear inverse
  problems}.
\bjtitle{SIAM Journal on Imaging Sciences}
\bvolume{2}(\bissue{1}),
\bfpage{183}--\blpage{202}
(\byear{2009})
\end{barticle}
\endbibitem

\bibitem[\protect\citeauthoryear{Dubhashi and Panconesi}{2009}]{dubhashi2009}
\begin{bbook}
\bauthor{\bsnm{Dubhashi}, \binits{D.P.}},
\bauthor{\bsnm{Panconesi}, \binits{A.}}:
\bbtitle{Concentration of Measure for the Analysis of Randomized Algorithms}.
\bpublisher{Cambridge University Press},
\blocation{Cambridge}
(\byear{2009})
\end{bbook}
\endbibitem

\end{thebibliography}

\end{document}